\documentclass[12pt,a4paper,draft]{article}
\usepackage[catalan,spanish,english]{babel}
\usepackage[psamsfonts]{amssymb}
\usepackage{cmmib57}
\usepackage{exscale}
\usepackage{amsmath}
\usepackage{amscd}
\usepackage{latexsym}
\usepackage{graphicx}
\newtheorem{thm}{Theorem}[section]
\newtheorem{cor}[thm]{Corollary}
\newtheorem{lem}[thm]{Lemma}
\newtheorem{prop}[thm]{Proposition}

\newtheorem{rem}{Remark}[section]
\newcommand{\qed}{\quad{$\blacksquare$}}

\newcommand{\Ind}{\;{\rm l}\hskip -0.23truecm 1}
\begin{document}
\begin{titlepage}
\null \vspace{2cm}
\begin{center}
{\Large\bf A fractional Poisson equation:\\[2mm]
existence, regularity and approximations \\[2mm]
of the solution \\[2mm]}
\medskip

by\\
\vspace{7mm}

\begin{tabular}{l@{\hspace{10mm}}l@{\hspace{10mm}}l}
{\sc Marta Sanz-Sol\'e}$\,^{(\ast)}$ &and &{\sc Iv\'an Torrecilla}$\,^{(\ast)}$\\
{\small marta.sanz@ub.edu }         &&{\small itorrecillatarantino@gmail.com}\\
{\small http://www.mat.ub.es/$\sim$sanz}&&\\
\end{tabular}
\begin{center}
{\small Facultat de Matem\`atiques}\\
{\small Universitat de Barcelona } \\
{\small Gran Via de les Corts Catalanes, 585} \\
{\small E-08007 Barcelona, Spain} \\
\end{center}
\end{center}

\vspace{1.5cm}

\noindent{\bf Abstract:} We consider a stochastic boundary value elliptic problem on a bounded domain
$D\subset \mathbb{R}^k$, driven by a fractional Brownian field with Hurst parameter $H=(H_1,\ldots,H_k)\in[\frac{1}{2},1[^k$.
First, we define the stochastic convolution derived from the Green kernel and prove some properties. Using monotonicity methods,
we prove existence and uniqueness of solution along with regularity of the sample paths. Finally, we propose a sequence of lattice approximations and
prove its convergence to the solution of the SPDE  at a given rate.
\medskip



\vspace{2 cm}

\noindent


{\begin{itemize} \item[$^{(\ast)}$] Supported by the grant MTM
2006-01351 from the \textit{Direcci\'on General de
Investigaci\'on, Ministerio de Ciencia e Innovaci\'on, Spain}
\end{itemize}}
\end{titlepage}

\newpage


\section{Introduction }
\label{s1}

This article deals with a stochastic Poisson equation on a bounded domain $D\subset \mathbb{R}^k$,
with arbitrary dimension $k\ge 1$, driven by a fractional Brownian field
$B^H$, with $H=(H_1,\ldots,H_k)\in[\frac{1}{2}, 1[^k$. We prove a theorem on existence
and uniqueness of solution, we study the properties of its sample paths and finally, we give a numerical scheme based on
lattice approximations and prove its convergence on a functional space with some explicit rate.

The equation is described as follows:
\begin{equation}
\label{1.1}
\left\{\renewcommand{\arraystretch}{1.8}
\begin{array}{rl}
\Delta u(x)-f(u(x))&=\ g(x)+\dot{B}^H(x),\ x\in D,\\
u(x)&=\ 0,\  x\in \partial D.
\end{array}
\right.
\end{equation}
We assume that $f$ has a decomposition $f=f_1+f_2$,
with $f_1,f_2:\mathbb{R}\to\mathbb{R}$ satisfying
\begin{description}
\item{(f1)} $f_1$ is continuous, non-decreasing and
$\sup_{x\in\mathbb{R}}|f_1(x)|\leq M$,
\item{(f2)} $f_2$ is Lipschitz with {\sl small} Lipschitz constant
$L$.
\end{description}
The function $g: D\to \mathbb{R}$ is measurable and satisfies some
integrability conditions. The stochastic
character of the equation comes from $\dot{B}^H(x)$, which denotes the
formal derivative of a fractional Brownian field.

We give a rigorous meaning to (\ref{1.1}) by means of a {\it mild} formulation,
as it is pretty usual in the SPDEs literature.
For this, we recall that if $k\ge 2$,  the Green function of the deterministic Poisson equation
on a bounded domain $D$ with smooth boundary is given by
\begin{equation}
\label{1.2} G_D^k(x,y)= G^k(x,y)-\mathbb{E}_x(G^k(B_\tau,y)),\ x,y\in D,
\end{equation}
with
\[
G^k(x,y)=C_k\left\{\renewcommand{\arraystretch}{1.8}
\begin{array}
[c]{ll}
\log |x-y|, & k=2,\\
|x-y|^{2-k}, & k\geq 3.
\end{array}
\right.
\]
Here $C_2=\frac{1}{2\pi}$, $C_k=\frac{1}{k(2-k)\omega_k}$ for $k\ge 3$,  where $\omega_k$ denotes the volume of the unit ball in $\mathbb{R}^k$, and $B_\tau$ is the random variable obtained by stopping a $k$-dimensional Brownian
motion starting at $x$ at its first exit time of
$D$ (see for instance \cite{GT97}  and also \cite{Do84}). For $k=1$, $G^{1}(x,y)= C_1|x-y|$ (see for instance \cite{K}, pg. 16). The expression for  $G^1_D(x,y)$
depends of the domain $D$. If $D=]0,1[$ then  $G_D^1(x,y)= (x\wedge y)-xy$, where ``$\wedge"$ denotes the infimum (see \cite{Do84} pg. 258).

By a solution to (\ref{1.1}) we mean a stochastic process
$u=\left\{u(x), x\in D\right\}$ satisfying
\begin{equation}
\label{1.3} u(x)=\int_D G_D^k(x,y)f(u(y))dy+\int_D G_D^k(x,y)g(y)dy+\int_D G_D^k(x,y)dB^H(y).
\end{equation}
For a similar SPDE in dimensions $k=1,2,3$, driven by a standard Wiener field $W$, different problems have been studied so far.
For instance, existence and uniqueness of solution has been proved in
\cite{BuPa90} using the classical theory of equations defined by monotone operators (see \cite{Ze90});
the Markov field property has been investigated in \cite{DM92} and \cite{DN94}, and numerical approximations
have been given in \cite{GM06}. Let us remark that for $k\le3$ the stochastic convolution
$\int_D G_D^k(x,y) dW(y)$ is well defined as a Wiener integral, because $G_D^k(x,\cdot)$
is square integrable. For $k\ge 4$, this property is not true anymore.

For $k\ge 4$, a SPDE of the same type
than (\ref{1.3}), driven by a Gaussian stationary process $F$ with an absolutely continuous covariance
measure, but possibly having singularities, has been studied in \cite{MSS06}, extending
the results of \cite{GM06} to higher dimensions. The authors combine conditions on
deterministic functions and covariance densities derived from Young's type inequalities
and provide a definition of an integral with respect to the random field $F$, and thereby a suitable meaning of the
stochastic convolution $\int_D G_D^k(x,y) dF(y)$. As regards the approximation scheme, the approach for $k\le 3$
based on a Fourier series expansion of $G_D^k(x,y)$ is not appropriate. In fact, as has been already mentioned,
$G_D^k(x,\cdot)$ is not square integrable and therefore this function and its Fourier series may not coincide.
 Instead, a more sophisticated procedure
involving a smoothing of the Green function combined with its Fourier series expansion has been considered.

The main reason for considering correlated noises in high dimensions is to compensate the irregularity of the Green function.
However, it may also be a reasonable choice when modeling phenomena where the stochastic input shows some dependence.

With the increasing attention devoted to fractional Brownian motion in recent years, the study of different type of
problems on SPDEs driven by fractional noise is being more present in the mathematical literature. We refer the
reader to \cite{SSV07} for an extensive list of references on the subject, including some motivating aspects from
other disciplines. Generically, these SPDEs are appropriate to model phenomena showing up either persistence
(for example, in hydrology or finance) or intermittency (as turbulence).
At the best of our knowledge, developments on this topic refer so far mainly to parabolic and hyperbolic
SPDEs, the elliptic case being less explored.

A  particular version of Equation (\ref{1.3}) with null functions $f$ and $g$ appears in \cite{HOZ00} (for $H_i\in]\frac{1}{2},1[^k$) and in \cite{HOZ04}
(for $H_i\in]0,1[^k$). In both references, the authors apply white noise analysis to give a meaning to the solution
$u(x)=\int_D G_D^k(x,y)dB^H(y)$ in the sense of distributions. Conditions on $H_i$ ensuring the existence of an $L^2(D)$-valued solution are given.
In comparison with these references, our analysis of (\ref{1.3}) allows a monotone nonlinearity $f(u)$ and a free term given by $g$,
as in \cite{GM06} and \cite{MSS06}.

Very recently, in \cite{LR08}, a new and very promising approach to stochastic Poisson equations based on Wiener chaos solutions on weighted spaces has been given.
Using Malliavin calculus and a formulation of the stochastic convolution by means of the Wick product, the authors prove existence and uniqueness of solution for elliptic equations
allowing some type of nonlinearities in the noisy term, modelled by a general spatial Gaussian process. Furthermore, in \cite{KRW09} numerical approximations based on a finite element
procedure are provided.

The content of the paper is as follows. In Section \ref{s2}, we
give some preliminaries on the fractional Brownian field $B^H$,
when $H=(H_1,\ldots,H_k)\in[\frac{1}{2}, 1[^k$. We combine ideas
from \cite{HOZ00} and \cite{HOZ04} (see also \cite{EH03}) with
some results from \cite{PT01} and \cite{Nu06} to give a moving
average representation of $B^H$ in terms of an standard Wiener
field. Then, we identify a suitable $L^p$-space with mixed norm of
deterministic functions which can be integrated against $B^H$.
These spaces are related with the reproducing kernel Hilbert space
of $B^H$ by means of Hardy-Sobolev's inequality. Section \ref{s3}
is devoted to the study of the stochastic convolution of the
Poisson kernel (\ref{1.2}). We give a sufficient condition on the
Hurst parameter $H$ ensuring the integrability of the Poisson
kernel with respect to $B^H$, according to the result proved in
Section \ref{s1}. We also give some probabilistic properties of
the stochastic convolution and prove the H\"older continuity of
its sample paths (Theorem \ref{t2.1}). These ingredients, though of own interest, are meant to
provide a rigorous meaning to Equation (\ref{1.3}) in any
dimension $k\ge 1$.
In Section \ref{s4}, we give a theorem on existence and
uniqueness of solution of Equation (\ref{1.3}) on the space of continuous functions vanishing at the boundary $\partial D$; we also prove H\"older continuity of the sample paths
of the solution. Finally, Section \ref{s5} is devoted to  numerical approximations of (\ref{1.3}). We consider the domain  $D=]0,1[^k$ and use the approach of \cite{GM06} for $k\le 3$ and
that of \cite{MSS06} when $k\ge 4$.
With an appropriate  choice of the functional spaces, we give the rate of convergence. For $k\le 3$ we find the same rate as for the Brownian case, while in dimensions $k\ge 4$, the rate depends on the regularity of the noise, as may be expected.

Throughout the paper, we shall denote by $c_H$ any positive constant depending on the Hurst parameter $H=(H_1,\ldots H_k)$, $k\ge 1$, independently of
its particular value, and by $C$ any positive, finite constant.


\section{Preliminaries}
\label{s2}
Let $H=(H_1,\ldots,H_k)\in ]0,1[^k$. A fractional Brownian field (fBf) on $\mathbb{R}^k$ with Hurst parameter $H$
is a Gaussian stochastic processes
$B^H=\{B^H(x),\;x\in\mathbb{R}^k\}$ with zero mean and  covariance function given by
\begin{equation}
\label{2.1}
R_H(x,y)=\mathbb{E}\left(B^H(x)B^H(y)\right)=\prod_{i=1}^kR_{H_i}(x_i,y_i),
\end{equation}
where
\begin{equation*}
R_{H_i}(x_i,y_i)=\frac{1}{2^k}\left(|y_i|^{2H_i}+|x_i|^{2H_i}-|x_i-y_i|^{2H_i}\right).
\end{equation*}
Such a process has been introduced and considered in relation with
different problems in \cite{Hu00}, \cite{Hu01}, \cite{HOZ00},
\cite{HOZ04}.

In this article, we consider values of the Hurst parameter
$H=(H_1,\ldots,H_k)\in [\frac{1}{2},1[^k$.
Our first goal is to define a stochastic convolution for the Poisson kernel with respect to $B^H$. For this, we shall identify a
suitable class of deterministic functions $f:\mathbb{R}^k\to \mathbb{R}$ for which
\begin{equation*}
\mathcal{I}^H(f)= \int_{\mathbb{R}^k} f(x)B^H(dx)
\end{equation*}
is a well defined  random variable.
As in the one parameter case, it will be useful to have a moving average type representation of the
process $B^H$ in terms of a standard Brownian field on $\mathbb{R}^k$. We shall prove such a representation
owing ideas from  \cite{HOZ04} but considering
the framework  of \cite{PT01} and \cite{Nu06} (see also Lemma 1.20.10 in \cite{mishura}).

We start by introducing some notation. On $\mathbb{R}^k$  we consider the usual partial order defined coordinatewise
and denote by $x=(x_1,\ldots,x_k)$ a generic element in this space.
For $x, y\in \mathbb{R}^k$ satisfying $x\le y$, we set
$\Ind_{[x,y]}(\eta)=\prod_{i=1}^k\Ind_{[x_i,y_i]}(\eta_i)$.
We shall denote by $\mathcal{E}$ the set of {\it elementary functions} on $\mathbb{R}^k$, that is
functions of the form
\begin{equation*}
\varphi(\eta)= \sum_{l=1}^{l_0} \varphi_l \Ind_{[x^l,y^l]}(\eta),
\end{equation*}
with $\varphi_l\in\mathbb{R}$ and disjoint rectangles $[x^l,y^l]$, $l=1,\ldots,l_0$. For $\varphi\in\mathcal{E}$,
we define
\begin{equation*}
\mathcal{I}(\varphi)= \sum_{l=1}^{l_0} \varphi_l B^H\left([x^l,y^l]\right),
\end{equation*}
where $B^H\left([x^l,y^l]\right)$ denotes the increment of $B^H$ on the rectangle $[x^l,y^l]$
in the sense of $k$--dimensional distribution functions.
For any $x\in\mathbb{R}^k$, we set $\Ind_{[0,x]}(\eta)=\prod_{i=1}^k\Ind_{[0,x_i]}(\eta_i)$,
where
\begin{equation*}
\Ind_{[0,x_i]}(\eta_i)=\left\{
\begin{array}{ll}
1, &  \eta_i\in[0,x_i], \\
-1, & \eta_i\in[x_i,0],\\
 0, &\text{otherwise}.
\end{array}
\right.
\end{equation*}
Then on $\mathcal{E}$, we introduce the inner product $\langle\cdot,\cdot\rangle_{\mathcal{H}^H}$ derived from the covariance structure of
$B^H$ given by
\begin{equation}
\label{2.2} \langle \Ind_{[0,x]}, \Ind_{[0,y]}
\rangle_{\mathcal{H}^H}=R_H(x,y).
\end{equation}
The closure of $\mathcal{E}$ with respect to the norm
$\Vert\cdot\Vert_{\mathcal{H}^H}$ will be denoted by $\mathcal{H}^H$.

For each $i=1,\ldots,k$,
we define the linear operator acting on functions $\varphi\in\mathcal{E}$ as follows:
\begin{equation*}
\label{2.3}
\left(K_{H_i}^\ast\varphi\right)(\eta)=\left\{\renewcommand{\arraystretch}{1.8}
\begin{array}{ll}
\varphi(\eta), &  H_i=\frac{1}{2},\\
c_{H_i}\int_{\mathbb{R}}\varphi(\eta_1,\ldots,\eta_{i-1},u_i,\eta_{i+1},\ldots,\eta_k)(u_i-\eta_i)_+^{H_i-\frac{3}{2}}du_i,
& H_i\in]\frac{1}{2},1[,
\end{array}
\right.
\end{equation*}
 where $c_{H_i}$ are constants depending only on $H_i$ (see (2.9) in \cite{HOZ04} for its explicit value).

For its further use, we introduce the sets $C^{(>)}= \{i=1,\ldots,k: H_i>\frac{1}{2}\}$,
$C^{(=)}= \{i=1,\ldots,k: H_i=\frac{1}{2}\}$, and denote by $c^{(>)}$, $c^{(=)}$, the cardinals
of $C^{(>)}$ and $C^{(=)}$, respectively. We notice that
\begin{equation}
\label{compact}
\left(K_{H_i}^\ast\varphi\right)(\eta)=\int_{\mathbb{R}}\varphi(\eta_1,\ldots,\eta_{i-1},u_i,\eta_{i+1},\ldots,\eta_k) \mu_i^\eta(du_i),
\end{equation}
where
\begin{equation*}
 \mu_i^\eta(du_i)=\left\{\renewcommand{\arraystretch}{1.8}
\begin{array}{ll}
c_{H_i}(u_i-\eta_i)_+^{H_i-\frac{3}{2}} du_i, &  i\in C^{(>)},\\
\delta_{\eta_i}(du_i), &  i\in C^{(=)},
  \end{array}
  \right.
  \end{equation*}
and $\delta_{\eta_i}$ denotes the Dirac measure at $\eta_i$.

By iteration, for $k\ge 2$ we define
\begin{equation}
\label{2.4}
\left(K_H^{\ast,(k)}\varphi\right)
=\left(K_{H_k}^\ast\left(K_{H_{k-1}}^{\ast}\cdots\left(K_{H_1}^\ast\varphi\right)\cdots\right)\right).
\end{equation}
Notice that for a function $\varphi=\otimes_{i=1}^k\varphi_i$, 
\begin{equation}
\label{product}
\left(K_H^{\ast,(k)}\varphi\right)(\eta) = \prod_{i=1}^k \left(K_{H_i}^\ast \varphi_i\right)(\eta_i).
\end{equation}

Consider the one-dimensional case ($k=1$). The Mandelbrot-van Ness representation establishes that
for any $H\in]0,1[$,
\begin{equation*}
B^H(x)= \int_{\mathbb{R}} K_H(x,y) W(dy),
\end{equation*}
where
\begin{equation}
\label{2.6}
K_{H}(x,y)= c_{H}\left\{\left(x-y\right)^{H-\frac{1}{2}}_+-(-y)^{H-\frac{1}{2}}_+\right\},
\end{equation}
and $W$ is a standard Brownian motion on the real line (see \cite{ST94}, pg. 320). Thus,
\begin{equation*}
\int_\mathbb{R}K_{H}(x,v)K_{H}(y,v)dv=R_{H}(x,y).
\end{equation*}
A simple computation yields
$\left(K_{H}^\ast\;\Ind_{[0,x]}\right)(y)=K_{H}(x,y)$, for $H>\frac{1}{2}$.
Hence,
\begin{equation*}
\mathcal{I}\left(\Ind_{[0,x]}\right):=B^H(x) = \int_{\mathbb{R}}\left(K^\ast_H\Ind_{[0,x]}\right)(y) W(dy).
\end{equation*}
With the definition (\ref{compact}) this representation also holds for  $H=\frac{1}{2}$.
Thus, it is clear that the mapping
$\Ind_{[0,x]}\mapsto \int_{\mathbb{R}} K_H^{\ast}\left(\Ind_{[0,x]}\right)(y) W(dy)$
is an isometry between $(\mathcal{E}, \langle\cdot,\cdot\rangle_{\mathcal{H}^H})$ and $L^2(\Omega)$.
That means, the operator $K_H^\ast$ is an isometry  between $(\mathcal{E}, \langle\cdot,\cdot\rangle_{\mathcal{H}^H})$
and $L^2(\mathbb{R})$.
Hence, $K_H^{\ast}$ can be extended to the Hilbert space $\mathcal{H}^H$. Otherwise stated,
we can extend the definition of $\mathcal{I}$ from $\mathcal{E}$ to $\mathcal{H}^H$, and therefore define
$ \mathcal{I}(\varphi)= \int_{\mathbb{R}}(K_H^{\ast}\varphi)(y)W(dy)$
as an $L^2(\Omega)$--valued random variable.
\smallskip

In the multidimensional case we prove similar results. We first recall a definition and introduce some notation.
A stochastic process $\{W(x), x\in \mathbb{R}^k\}$ is termed a
{\it standard Wiener field} on $\mathbb{R}^k$ if it is Gaussian, with mean zero
and covariance given by
$\mathbb{E}\left(W(x)W(y)\right) = x\wedge y$,
where $x\wedge y:=\prod_{i=1}^k (x_i\wedge y_i)$, and
\begin{equation*}
x_i\wedge y_i=\left\{
\begin{array}{ll}
x_i\wedge y_i, & x_i,y_i>0, \\
(-x_i)\wedge (-y_i), & x_i,y_i<0,\\
0, &\text{otherwise}.
\end{array}
\right.
\end{equation*}
For a function $\varphi:\mathbb{R}^k\to \mathbb{R}$, we define
$\tilde\varphi(u_1,\ldots,u_k;v_1,\ldots,v_k)=\varphi(w_1,\ldots,w_k)$, with $w_i=u_i$, if $i\in C^{(=)}$ and $w_i=v_i$ if $i\in C^{(>)}$. Then we denote by $\left|\mathcal{H}^H\right|$
the set of functions $\varphi:\mathbb{R}^k\to \mathbb{R}$ such that
\begin{align*}
& \int_{\mathbb{R}^{c^{(=)}}} \prod_{i\in C^{(=)}} du_i  \int_{\mathbb{R}^{2c^{(>)}}} \prod_{i\in C^{(>)}} \left(du_i dv_i H_i(2H_i-1)|u_i-v_i|^{2H_i-2}\right)\\
&\quad\times|\varphi(u_1,\ldots,u_k)| |\tilde\varphi(u_1,\ldots,u_k;v_1,\ldots,v_k)|<+\infty.
\end{align*}

\begin{prop}
\label{p2.1} Set $K_H(x,y)= \prod_{i=1}^k K_{H_i}(x_i,y_i)$, with
$K_{H_i}$, $i=1,\ldots,k$, defined in (\ref{2.6}). Then:
\begin{enumerate}
\item There exists a
standard Wiener field $W$ on $\mathbb{R}^k$ such that
\begin{equation}
\label{2.7}
B^H(x)= \int_{\mathbb{R}^k} K_H(x,y) W(dy) =  \int_{\mathbb{R}^k} \left(K^{\ast,(k)}_H\Ind_{[0,x]}\right)(y)W(dy).
\end{equation}
\item
For any
$\varphi \in \mathcal{H}^H$
\begin{equation}
\label{int}
\mathcal{I}(\varphi)=\int_{\mathbb{R}^k}\left(K_H^{\ast,(k)}\varphi\right)(y)dW(y)
\end{equation}
defines a random variable in $L^2(\Omega)$.
\item  For any $\varphi_1, \varphi_2\in \left|\mathcal{H}^H\right|$ the following isometry formula holds:
\begin{align}
\label{isometry}
\mathbb{E}\left(\mathcal{I}(\varphi_1)\mathcal{I}(\varphi_2)\right)&=\int_{\mathbb{R}^{k}}\left(K_H^{\ast,(k)}\varphi_1\right)(y)
\left(K_H^{\ast,(k)}\varphi_2\right)(y)dy\nonumber\\
&= \int_{\mathbb{R}^{c^{(=)}}} \prod_{i\in C^{(=)}} du_i  \int_{\mathbb{R}^{2c^{(>)}}} \prod_{i\in C^{(>)}} \left(du_i dv_i H_i(2H_i-1)|u_i-v_i|^{2H_i-2}\right)\nonumber\\
&\quad \times \varphi_1(u_1,\ldots,u_k)\tilde\varphi_2(u_1,\ldots,u_k;v_1,\ldots,v_k).
\end{align}

\noindent In particular, if $H_i\in]\frac{1}{2},1[$ for any $i=1,\ldots,k$,
\begin{equation}
\label{isometrybis}
\mathbb{E}\left(\mathcal{I}(\varphi_1)\mathcal{I}(\varphi_2)\right)= \int_{\mathbb{R}^{2k}} \varphi_1(u)\varphi_2(v) \prod_{i=1}^k\left(H_i(2H_i-1) |u_i-v_i|^{2H_i-2}\right) du dv.
\end{equation}
\end{enumerate}
\end{prop}

\noindent{\it Proof}: To prove the  existence of $W$, we follow the arguments of \cite{Nu06}, pg. 279, for $k=1$, which extend easily to any $k\ge 1$, as follows.

For each $i=1,\ldots,k$, such that $H_i>\frac{1}{2}$, the action of the
kernel $K_{H_i}^\ast$  on elementary functions can be expressed in terms of a
fractional integral. More precisely,
if $\varphi\in\mathcal{E}$,
$K_{H_i}^\ast\varphi=c_{H_i}I_-^{H_i-\frac{1}{2}}\varphi$, where $I_-^{H_i-\frac{1}{2}}\varphi(x)= \frac{1}{\Gamma(H_i-\frac{1}{2})}\int_x^{\infty}(y-x)^{H_i-\frac{3}{2}}\varphi(y) dy$
(see section 3.2 in \cite{PT01}).
Then, for $\psi$ in the image of  $K_{H_i}^\ast$ and by considering the fractional derivative defined by
$D_-^{H_i-\frac{1}{2}}\psi(x) = \frac{1}{\Gamma(\frac{3}{2}-H_i)}\int_0^\infty \left(\psi(x)-\psi(x+y)\right) y^{-\frac{1}{2}-H_i} dy$
we define
$Q_{H_i}^\ast\psi=c_{H_i}^{-1}D_-^{H_i-\frac{1}{2}}\psi$,
which by  the rules of fractional calculus is seen to be  the
inverse operator of $K_{H_i}^\ast$.
For  $H_i=\frac{1}{2}$,  $Q_{H_i}^\ast$ is defined to be the identity operator.

Set $W(x)= B^H \left(\prod_{i=1}^k Q_{H_i}^\ast \Ind_{[0,x_i]}\right)$.
The process
$W=\left\{W(x),\;x\in\mathbb{R}^k\right\}$ is a standard Wiener field on $\mathbb{R}^k$, and the stochastic field $B^H$ has the
integral representation
\begin{equation}
\label{2.14} B^H(x)=\int_{\mathbb{R}^k}K_H(x,y)dW(y).
\end{equation}
Indeed, set $Q_H^{\ast,(k)}\Ind_{[0,x]}= \prod_{i=1}^k Q_{H_i}^\ast \Ind_{[0,x_i]}$. Using (\ref{product}), we have
$\left(K_H^{\ast,(k)}\circ Q_H^{\ast,(k)}\right)\Ind_{[0,x]}=\Ind_{[0,x]}$. Thus,
for any $x,y\in\mathbb{R}^k$,
\begin{align*}
\mathbb{E}\left(W(x)W(y)\right)
&=\left \langle
Q_H^{\ast,(k)}\Ind_{[0,x]},Q_H^{\ast,(k)}\Ind_{[0,y]}\right
\rangle_{\mathcal{H}^H}\notag\\
&=\left \langle \Ind_{[0,x]}, \Ind_{[0,y]}\right
\rangle_{L^2(\mathbb{R}^k)}=x\wedge y.
\end{align*}
By construction, the operator $K_H^{\ast,(k)}$ is an isometry from $\mathcal{E}$ into $L^2(\mathbb{R}^k)$ that can be extended to the
Hilbert space $\mathcal{H}^H$. Therefore, one can define $\mathcal{I}(\varphi)$ for any $\varphi\in\mathcal{H}^H$ by means of (\ref{int}).
\smallskip

We now prove (\ref{isometry}).  By the very definition of
$K_H^{\ast,(k)}$ (see (\ref{compact})) and by applying
Fubini's theorem, we obtain
\begin{align*}
&\int_{\mathbb{R}^k}\left(K_H^{\ast,(k)}\varphi_1\right)(y)\left(K_H^{\ast,(k)}\varphi_2\right)(y)dy\\
&=\int_{\mathbb{R}^k} dy \left(\int_{\mathbb{R}^k} \varphi_1(u) \prod_{i\in C^{(=)}} \mu_i^y(du_i) \prod_{i\in C^{(>)}} \mu_i^y(du_i)\right)\\
&\quad \times \left(\int_{\mathbb{R}^k} \varphi_2(v) \prod_{i\in C^{(=)}} \mu_i^y(dv_i) \prod_{i\in C^{(>)}} \mu_i^y(dv_i)\right)\\
&=\int_{\mathbb{R}^{c^{(=)}}} \prod_{i\in C^{(=)}} du_i \\
&\quad \times  \int_{\mathbb{R}^{2c^{(>)}}}\prod_{i\in C^{(>)}} du_i dv_i\left( \int_{\mathbb{R}^{c^{(>)}}}\prod_{i\in C^{(>)}} dy_i c_{H_i}^2
(u_i-y_i)_+^{H_i-\frac{3}{2}} (v_i-y_i)_+^{H_i-\frac{3}{2}}\right)\\
&\quad \times \varphi_1(u_1,\ldots,u_k)\tilde\varphi_2(u_1,\ldots,u_k;v_1,\ldots,v_k).
\end{align*}
From this and the identity
\begin{equation*}
c_{H_i}^2 \int_{\mathbb{R}} dy_i (u_i-y_i)_+^{H_i-\frac{3}{2}}(v_i-y_i)_+^{H_i-\frac{3}{2}}
=H_i(2H_i-1)\vert u_i-v_i\vert^{2H_i-2}
\end{equation*}
(see \cite{GN96}, page 404),
(\ref{isometry}) follows.

Finally, if $C^{(>)}=\{1,\ldots,k\}$, (\ref{isometry}) reads (\ref{isometrybis}).
This ends the proof of the Proposition.


\hfill\qed

\medskip


It is well known that for real functions $\varphi$, $\psi$ and $H\in]\frac{1}{2},1[$,
\begin{equation}
\label{2.15}
\int_{\mathbb{R}}\int_{\mathbb{R}}|\varphi(\eta)||\psi(\theta)| |\eta-\theta|^{2H-2}d\eta
d\theta\leq b_{H}
\|\varphi\|_{L^\frac{1}{H}(\mathbb{R})}\|\psi\|_{L^\frac{1}{H}(\mathbb{R})},
\end{equation}
with a positive constant $b_{H}$. Indeed, this follows from Hardy-Littlewood-Sobolev's inequality
(see inequality $(1)$, page 321, in \cite{BePa61}).
By considering functions $\varphi_1, \varphi_2:\mathbb{R}^k\to \mathbb{R}$ and
applying recursively this inequality for indices $i\in C^{(>)}$, and Schwarz inequality for $i\in C^{(=)}$, we obtain
\begin{align}
\label{2.16}
& \int_{\mathbb{R}^{c^{(=)}}} \prod_{i\in C^{(=)}} du_i  \int_{\mathbb{R}^{2c^{(>)}}} \prod_{i\in C^{(>)}} \left(du_i dv_i H_i(2H_i-1)|u_i-v_i|^{2H_i-2}\right)\nonumber\\
&\quad\times|\varphi_1(u_1,\ldots,u_k)||\tilde\varphi_2(u_1,\ldots,u_k;v_1,\ldots,v_k)|\nonumber\\
&\le c_{H_1,\ldots,H_k} \Vert
\varphi_1\Vert_{L^{\frac{1}{H_1},\ldots,\frac{1}{H_k}}(\mathbb{R}^k)}\Vert
\varphi_2\Vert_{L^{\frac{1}{H_1},\ldots,\frac{1}{H_k}}(\mathbb{R}^k)},
\end{align}
 where for $p_i\in[1,\infty]$, $i=1,\ldots,k$,
\begin{equation*}
\|h\|_{L^{p_1,\ldots,p_k}(\mathbb{R}^k)}
=\Big(\int_{\mathbb{R}}\Big(\int_{\mathbb{R}}\cdots\Big(\int_{\mathbb{R}}|h(\eta_1,\ldots,\eta_k)|^{p_1}d\eta_1\Big)^{\frac{p_2}{p_1}}\cdots\Big)^{\frac{p_k}{p_{k-1}}}d\eta_k\Big)^{\frac{1}{p_k}}.
\end{equation*}
Let us denote by $L^{p_1,\ldots,p_k}(\mathbb{R}^k)$
the space of measurable functions $h$ defined on $\mathbb{R}^k$ with $\|h\|_{L^{p_1,\ldots,p_k}(\mathbb{R}^k)}<\infty$.
Such  spaces are termed $L^p$ {\it spaces with mixed norm}. For details we refer the reader to \cite{BePa61}
and also \cite{AF03}. In particular, if $H_i\in]\frac{1}{2},1[$ for any $i=1,\ldots,k$, a proof of (\ref{2.16}) is given in page 322 of \cite{BePa61},
but it is easy to extend the result allowing the value $H_i=\frac{1}{2}$ for some indices $i$.
The preceding discussion yields
\begin{equation}
\label{2.18}
L^{\frac{1}{H_1},\ldots,\frac{1}{H_k}}(\mathbb{R}^k)\subset \left\vert\mathcal{H}^H\right\vert\subset \mathcal{H}^H.
\end{equation}

For its further use, we remark that for any $p\ge\sup_{i\in\{1,\ldots,k\}}p_i$, and every measurable function $h$ with bounded support
$\mathcal{O}$ contained in $\mathbb{R}^k$,
\begin{equation}
\label{2.17}
\|h\|_{L^{p_1,\ldots,p_k}(\mathcal{O})}\le C \Vert h\Vert_{L^p(\mathcal{O})},
\end{equation}
with a constant $C$ depending only on $\mathcal{O}$. Indeed, this follows by applying recursively H\"older's inequality with $pp_k^{-1}, \ldots, pp_1^{-1}$.


\section{The fractional stochastic convolution of the\break Poisson kernel}
\label{s3}
In this section, we consider a bounded domain with $\mathcal{C}^\infty$ boundary, $D\subset \mathbb{R}^k$, if $k\ge 2$, and
$D=]0,1[$ if $k=1$. We consider the Green function
defined in (\ref{1.2}) and $G^1_D(x,y)= (x\wedge y)-xy$, respectively.
The purpose is to define the stochastic
convolution $\int_D G_D^k(x,y) dB^H(y)$ with respect to the fractional Brownian field with
parameters $H_i\in[\frac{1}{2},1[$, $i=1,\ldots,k$, introduced in
the preceding section, and to study its sample paths.

Throughout the section we shall make use of the following remark:
Let $k\ge 2$ and assume that for some norm $\Vert\cdot\Vert$ defined on a space of functions on $\mathbb{R}^k$, we have $\sup_{x\in D}\Vert
G^k(x,\cdot)\Vert <\infty$. Then, $\sup_{x\in D}\Vert\mathbb{E}_x\left(G^k(B_\tau,\cdot)\right)\Vert <\infty$, and
consequently, $\sup_{x\in D}\Vert G_D^k(x,\cdot)\Vert <\infty$. Indeed,
since $\mathbb{E}_x$ is a convex operator, denoting by $P_x$
the probability law of $B_\tau$, we obtain
\begin{equation*}
\left\Vert \mathbb{E}_x\left(G^k(B_\tau,\cdot)\right)\right\Vert \le
\mathbb{E}_x\left\Vert G^k(B_\tau,\cdot)\right\Vert =\int_{\mathbb{R}^k}
P_x(dz) \left\Vert G^k(z,\cdot)\right\Vert \le \sup_{z\in D}\left\Vert G^k(z,\cdot)\right\Vert.
\end{equation*}

We are interested in the integrability properties of $G^k_D$.
To start with, let us state a result that for dimensions $k\ge 3$ is Lemma 2 in \cite{MSS06}.
The extension to $k=1, 2$ is trivial.

\begin{lem}
\label{l3.1} For any $p\in[1,\frac{k}{(k\vee 2)-2}[$, there exists a
positive constant $\mathcal{K}_1$ depending on $p$ and $k$, such
that
\begin{equation}
\label{3.1}
\sup_{x\in D}\|G_D^k(x,\cdot)\|_{L^p(D)}\leq \mathcal{K}_1.
\end{equation}
\end{lem}

For the values $k=1,2$, (\ref{3.1}) holds for any $p\in[1,\infty[$. Therefore by virtue of (\ref{2.17})
we can choose $p_0\ge\sup_{i\in\{1,\ldots,k\}}\frac{1}{H_i}$ such that
\begin{equation}
\label{3.2}
\sup_{x\in D}\Vert G^k_D(x,\cdot)\Vert_{L^{\frac{1}{H_1},\ldots\frac{1}{H_k}}(D)}\le C \sup_{x\in D}\Vert G^k_D(x,\cdot)\Vert_{L^{p_0}(D)}<\infty.
\end{equation}
Consider now the case $k=3$. Property (\ref{3.1}) holds for any $p\in[1,3[$. Hence there exist
$p_0\in [\sup_{i\in\{1,\ldots,k\}}\frac{1}{H_i}, 3[$ such that (\ref{3.2}) holds.
Then, accordingly with the results stated in the preceding section, for $k=1,2,3$, $\mathcal{I}(G^k_D(x,\cdot))$ is a well-defined  random variable in $L^2(\Omega)$ for any $x\in D$.

A similar conclusion holds true for $k=4$ under the additional assumption $H_i\in]\frac{1}{2}, 1[$, that is, excluding the possibility of
having a standard Brownian motion in some of the components of $B^H$. By similar arguments, the existence of $\mathcal{I}(G^k_D(x,\cdot))$
for $k\ge 4$ is ensured by the stronger hypothesis $H_i\in]\frac{k-2}{k},1[$, for any $i=1,\ldots,k$.
As we show in the next Lemma, one can relax this assumption by working with $L^p$ spaces with mixed norm.

Throughout the section we suppose that $D\subset [-R,R]^k$ for some $R>0$. We shall use the following inequality
on Euclidean norms:
For any $\mu\ge 0$ and $\beta_i\ge 0$, $i=1,\ldots,k$, such that $\sum_{i=1}^k \beta_i=1$,
\begin{equation}
\label{ineq}
|x|^{-\mu}\le \prod_{i=1}^k |x_i|^{-\beta_i\mu}.
\end{equation}

\begin{lem}
\label{l3.2}
Let $k\ge 4$ and assume that $\sum_{i=1}^k H_i>k-2$. Then, there
exists a positive constant $\mathcal{K}_2$ depending on $H$, $D$ and $k$,
such that
\begin{equation}
\label{3.3}
\sup_{x\in D}\|G_D^k(x,\cdot)\|_{L^{\frac{1}{H_1},\cdots,\frac{1}{H_k}}(D)}\leq \mathcal{K}_2.
\end{equation}
\end{lem}

\noindent {\it Proof}:
 By applying the inequality (\ref{ineq}) with $\mu:=2-k$ and $\beta_i:=\frac{H_i}{\sum_{i=1}^k H_i}$, and the remark at the beginning of the section, we obtain
\begin{align*}
\|G_D^k(x,\cdot)\|_{L^{\frac{1}{H_1},\cdots,\frac{1}{H_k}}(D)}&\le
2\Big(\int_{-R}^R\cdots\Big(\int_{-R}^R|x-y|^{\frac{2-k}{H_1}}dy_1\Big)^{\frac{H_1}{H_2}}\cdots
dy_k\Big)^{H_k}\\
& \leq 2\prod_{i=1}^k
\Big(\int_{-R}^R|x_i-y_i|^{\frac{\beta_i(2-k)}{H_i}}dy_i\Big)^{H_i}.
\end{align*}
The supremum on $x\in D$ of the last term is finite if and only if
$\frac{\beta_i(2-k)}{H_i}> -1$. By the definition of $\beta_i$,
this condition is equivalent to $\sum_{i=1}^k H_i>k-2$. \hfill\qed
\medskip

In the sequel we will assume the hypothesis:
\smallskip

\noindent {\bf(H)}\quad $H_i\in[\frac{1}{2},1[$, $i=1,\ldots,k$, and for dimensions $k\ge 4$, we suppose in addition that
$\sum_{i=1}^k H_i>k-2$.
\smallskip

Then,  for any dimension $k\ge 1$, Lemmas \ref{l3.1} and \ref{l3.2} yield
\begin{equation}
\label{3.3.1} \sup_{x\in
D}\|G_D^k(x,\cdot)\|_{L^{\frac{1}{H_1},\cdots,\frac{1}{H_k}}(D)}\le
\mathcal{K},
\end{equation}
with  $\mathcal{K}=\max( \mathcal{K}_1, \mathcal{K}_2)$.

This
proves the existence of a stochastic process
\begin{equation}
\label{3.4}
\mathcal{J}= \left\{\mathcal{J}(x) :=  \mathcal{I}\left(G_D^k(x,\cdot)\right)= \int_DG_D^k(x,y) B^H(dy), x\in D\right\},
\end{equation}
with values in $L^2(\Omega)$, satisfying
\begin{align}
\label{3.5} \sup_{x\in
D}\mathbb{E}\left|\int_{D}G_D^k(x,y)dB^H(y)\right|^2&
=\sup_{x\in D}\int_{D}\left(K_H^{\ast,(k)}G_D^k(x,\cdot)\right)^2(y)dy \nonumber\\
& \leq C_H\sup_{x\in
D}\|G_D^k(x,\cdot)\|^2_{L^{\frac{1}{H_1},\cdots,\frac{1}{H_k}}(D)}
\leq C_H{\mathcal{K}}^2< \infty.
\end{align}



Next we prove that the stochastic field $\mathcal{J}$ has a.s. H\"older
continuous sample paths.

\begin{thm} Under {\bf(${\rm \bf H}$)} it holds that
\label{t2.1}
\begin{equation}
\label{holder}
\left\Vert G_D^k(x,\cdot)-G_D^k(z,\cdot)\right\Vert_{L^{\frac{1}{H_1},\cdots,\frac{1}{H_k}}(D)}
\le C |x-z|^\lambda,
\end{equation}
for any $x,z\in D$, with
\begin{equation*}
\left\{\renewcommand{\arraystretch}{1.8}
\begin{array}{ll}
\lambda=1,&\text{ for
}k=1,\\
\lambda\in\left]0,\left(2-k+\sum_{i=1}^k H_i\right)\wedge 1\right[, &\text{
for }k \ge 2.
\end{array}
\right.
\end{equation*}
Therefore, the Gaussian random field $\mathcal{J}$
defined in $(\ref{3.4})$ satisfies
\begin{equation}
\label{3.7} \mathbb{E}\left(|\mathcal{J}(x)-\mathcal{J}(z)|^2\right)\leq
C|x-z|^{2\lambda},
\end{equation}
and a.s., the sample paths are H\"older continuous of order
$\gamma \in ]0,\lambda[$. Moreover, for any $p\in [1,\infty[$, $\mathcal{J}\in L^p\left(\Omega;L^\infty(D)\right)$.
\end{thm}

\noindent{\it Proof}: Fix $x,z\in D$. For $k=1$ and $D=]0,1[$, direct computations yield
\begin{equation*}
\left|G_D^k(x,y)-G_D^k(z,y)\right|\leq C|x-z|,
\end{equation*}
which implies (\ref{holder}).

Let $k\ge 2$ and set
\begin{align}
\label{3.7.1}
T^{(k)}&=\left\|G^k(x,\cdot)-G^k(z,\cdot)\right\|^2_{L^{\frac{1}{H_1},\cdots,\frac{1}{H_k}}(D)},\nonumber\\
\tilde T^{(k)}&=\left\|\mathbb{E}_x\left(G^k(B_\tau,\cdot)\right)-\mathbb{E}_z\left(G^k(B_\tau,\cdot)\right)\right\|^2_{L^{\frac{1}{H_1},\cdots,\frac{1}{H_k}}(D)}.
\end{align}
By the strong Markov property of Brownian motion,
$$\mathbb{E}_z\left(G^k(B_\tau,\cdot)\right)=\mathbb{E}_x\left(G^k(B_\tau-x+z,\cdot)\right).$$
Thus, similarly as in the remark at the beginning of the section,
\begin{equation*}
\tilde T^{(k)}\le \sup_{y\in D}\left\Vert
G^k(y,\cdot)-G^k(y-x+z,\cdot)\right\Vert^2_{L^{\frac{1}{H_1},\cdots,\frac{1}{H_k}}(D)}.
\end{equation*}
Therefore, we can concentrate on proving that
$T^{(k)}\leq C|x-z|^{2\lambda}$ for the values of $\lambda$ given in the
statement.

By the very definition of $G^k$,
\begin{equation*}
 \left|G^k(x,y)-G^k(z,y)\right|\le C \left\vert\int_{|z-y|}^{|x-y|} t^{1-k} dt\right\vert.
 \end{equation*}
Fix $\lambda\in]0,1]$ and consider the change of variable $t:= \left(\theta|x-y|^\lambda+(1-\theta)|z-y|^\lambda\right)^{\frac{1}{\lambda}}$,
$\theta\in[0,1]$. Then,
\begin{align}
\label{3.8.2} \left|G^k(x,y)-G^k(z,y)\right|&\leq C \left\vert \vert x-y\vert^\lambda -\vert z-y\vert^\lambda\right\vert \\
&\quad \times
\left|\int_0^1\left(\theta|x-y|^\lambda+(1-\theta)|z-y|^\lambda\right)^\frac{2-k-\lambda}{\lambda}d\theta\right|\nonumber\\
&\leq
C|x-z|^\lambda\left(|x-y|^{2-k-\lambda}+|z-y|^{2-k-\lambda}\right),
\end{align}
where in the last upper bound we have used the inequality $(a+b)^p\le a^p+b^p$, valid for any $a,b\ge 0$ and $0<p\le 1$.

We can now apply Minkowski's inequality for $L^p$ spaces with mixed norm (see \cite{BePa61}, page 302) to obtain
\begin{align*}
T^{(k)}&\leq C |x-z|^{2\lambda}\nonumber\\
&\times\Big(\Big\||x-\cdot|^{2-k-\lambda}\Big\|^2_{L^{\frac{1}{H_1},\cdots,\frac{1}{H_k}}(D)}+
\Big\||z-\cdot|^{2-k-\lambda}\Big\|^2_{L^{\frac{1}{H_1},\cdots,\frac{1}{H_k}}(D)}\Big).
\end{align*}
Suppose that $\lambda<2-k+\sum_{i=1}^kH_i$. Following the arguments of the proof of Lemma \ref{l3.2},
the inequality
(\ref{ineq}) applied to $\mu:=\lambda+k-2$ yields
\begin{equation*}
T^{(k)} \le C |x-z|^{2\lambda}.
\end{equation*}
This ends the proof of (\ref{holder}) and, by the isometry
property of the stochastic integral given in (\ref{3.5}), we
obtain (\ref{3.7}).

Since the process $\mathcal{J}$ is Gaussian, the statement about
regularity of its sample paths follows from Kolmogorov's
continuity criterion. More precisely, we get
\begin{equation*}
\sup_{x,z\in D, x\ne
z}\frac{\left\vert\mathcal{J}(x)-\mathcal{J}(z)\right\vert}{|x-z|^\gamma}:=
C(\omega)<+\infty,\; a.s.,
\end{equation*}
for any $\gamma\in]0,\lambda[$ (see \cite{R-Y}, Theorem 2.1, pg. 26).

From here, we easily obtain $\sup_{x\in D}\left\vert\mathcal{J}(x)\right\vert < + \infty$, a.s. This implies that $\mathcal{J}\in L^p(\Omega;L^\infty(D)$, by
Theorem 3.2 of \cite{adler}.
 \hfill \qed

 \begin{rem}
\label{r3.1}
Consider the following assumption which is stronger than {\bf(H)}:
\smallskip

\noindent{\bf(${\rm \bf H}^{\ast}$)}\quad
$H_i\in[\frac{1}{2},1[^k$, $i=1,\ldots,k$, and $\sum_{i=1}^kH_i>k-1$ for any $k\geq 2$,
\smallskip

\noindent
In this case, it can be proved that (\ref{holder}) holds true with $\lambda=1$ for any dimension $k\ge 1$. \end{rem}


\section{Existence and uniqueness of solution to the \break fractional Poisson equation}
\label{s4}

This section is devoted to establish  the existence and uniqueness
of solution to the equation (\ref{1.1}). This result will be
obtained by a pathwise argument; once it will be established, we
will prove some probabilistic properties of the solution. We
borrow the method of the proof from \cite{MSS06} (see also
\cite{BuPa90}, \cite{DM92} and \cite{GM06}), which follows the
classical monotonicity methods.
We shall denote by 
$\mathcal{C}(\bar D)$ the space of continuous functions on $\bar D$ and set
$\mathcal{S}=\left\{w: w \in \mathcal{C}(\bar D), w|_{\partial D}=0\right\}$.

For its further use we highlight some properties. The first one,
denoted by {\bf (M)} is a monotonicity property. The second one,
named  {\bf (P)}, has been proved in \cite{BuPa90} (Lemma $2.4$);
it is a consequence of the solvability of the Dirichlet problem on
$D$ and  Poincar\'e's inequality (see \cite{GT97} or \cite{AF03}).
They are formulated as follows:
\smallskip

{\sl {\bf (M)} $f$ is a function of the form $f=f_1+f_2$ with
$f_1,f_2:\mathbb{R} \to \mathbb{R}$, $f_1$ non-decreasing and
$f_2$ Lipschitz with Lipschitz constant $L$, if and only if for
every $u,v\in \mathbb{R}$,
\begin{equation}
\label{2.0} (u-v)(f(u)-f(v))\geq -L(u-v)^2.
\end{equation}
}
\smallskip

{\sl {\bf (P)} There exists a constant $a>0$ such that for any
$\varphi \in L^2(D)$,
\begin{equation}
\label{2.6.1} \int_D\left(\int_D
G_D^k(x,y)\varphi(y)dy\right)\varphi(x)dx\leq -a\int_D\left(\int_D
G_D^k(x,y)\varphi(y)dy\right)^2dx.
\end{equation}
}

We begin with the existence and uniqueness result.
\begin{thm}
\label{t4.1} Assume {\bf(H)} and
suppose that $g$ is a real function defined on $D$, $g\in L^{\frac{1}{1-H_1},\cdots,\frac{1}{1-H_k}}(D)$, $f=f_1+f_2$ 
 and satisfies $(f1)$ and $(f2)$ with a
Lipschitz constant  $L<a$. 
Then, there
exists a unique stochastic process solution to $(\ref{1.3})$ with sample paths in $\mathcal{S}$, a.s.
\end{thm}

\noindent{\it Proof}:
Consider the operator $\mathcal{T}:\mathcal{S} \longrightarrow
\mathcal{S}$, defined by
\begin{equation*}
\mathcal{T}(w)(x)=w(x)-\int_D G_D^k(x,y)f(w(y))dy.
\end{equation*}
H\"older's inequality for $L^p$ spaces with mixed norm and
(\ref{holder}) yield
\begin{align}
\label{4.1}
&\left|\int_D\left(G_D^k(x,y)-G_D^k(z,y)\right)g(y)dy\right|
\notag \\
&\leq\left\|G_D^k(x,\cdot)-G_D^k(z,\cdot)\right\|_{L^{\frac{1}{H_1},\cdots,\frac{1}{H_k}}(D)}
\left\|g\right\|_{L^{\frac{1}{1-H_1},\cdots,\frac{1}{1-H_k}}(D)}\notag\\
&\leq C|x-z|^\lambda,
\end{align}
with $\lambda>0$ given in Theorem $\ref{t2.1}$. With this and the continuity of the stochastic convolution established in Theorem $\ref{t2.1}$, we conclude
\begin{equation}
\label{defb}
b(x):=\int_D G_D^k(x,y)g(y)dy+\int_D G_D^k(x,y)dB^H(y)\in
\mathcal{S},\  a.s.
\end{equation}

We next show that the operator equation $\mathcal{T}w=b$ has
a unique solution for any $b\in\mathcal{S}$, or equivalently that
$\mathcal{T}$ is bijective. Uniqueness
guarantees the measurability of the process $\left\{w(x),\;x\in D\right\}$.

The one-to-one property of $\mathcal{T}$ follows by applying {\bf
(M)} and {\bf {(P)}} (see Theorem $2$ in \cite{MSS06}). Indeed,
let $u$, $v$ be such that $\mathcal{T}u=\mathcal{T}v$. Then
$u(x)-v(x)=\int_D G_D^k(x,y)(f(u(y))-f(v(y)))dy$. Multiplying both
sides of this equation by $f(u(x))-f(v(x))$, integrating over $D$
and applying property {\bf (P)}, we obtain
\begin{align*}
&\int_D\left(u(x)-v(x)\right)\left(f(u(x))-f(v(x))\right)dx\\
&\quad =\int_D\left(f(u(x))-f(v(x))\right)\left(\int_DG_D^k(x,y)\left(f(u(y))-f((y))\right) dy\right) dx\\
&\quad \le -a\int_D\left(\int_D G_D^k(x,y)\left(f(u(y))-f(v(y))\right)dy\right)^2 dx\\
&\quad = -a\int_D\left(u(x)-v(x)\right)^2 dx.
\end{align*}
Because of {\bf {(M)}},
\begin{equation*}
\int_D\left(u(x)-v(x)\right)\left(f(u(x))-f(v(x))\right)dx\ge -L\int_D\left(u(x)-v(x)\right)^2 dx.
\end{equation*}
Hence, $(a-L)\int_D\left(u(x)-v(x)\right)^2 dx \le 0$. Since $L<a$ and $u,v\in L^2(D)$, this implies $u(x)=v(x)$ for almost every $x\in D$.
\medskip

We now prove that $\mathcal{T}$ is onto. In the next argument, we fix a sample path of the process
$B^H$ on a set of probability one. 
\smallskip

\noindent{\sl Step 1: A solution for a regular problem.}

\noindent We follow the arguments of \cite{MSS06} Lemma 3 and the sequent discussion.  Let $b\in
\mathcal{S}$ and $b_n\in C_{c}^\infty(D)$, $n\geq 1$, be such that
$b_n\longrightarrow b$ in $L^2(D)$. Then, one can construct a
sequence of functions solving $\mathcal{T}u^{(n)}=b_n$ and
satisfying
\begin{equation}
\label{4.2} u^{(n)}(x)=\int_D
G_D^k(x,y)f(u^{(n)}(y))dy+b_n,\quad\text{for}\quad x\in D,\quad
u^{(n)}|_{\partial D}=0.
\end{equation}
By the properties {\bf (M)}, {\bf (P)} and since $u^{(n)}\in
L^2(D)$, one can prove that $\{u^{(n)},\;n\geq 1\}$ is a Cauchy
sequence in $L^2(D)$. Let $u$ denote the limit.
\smallskip

\noindent{\sl Step 2: $u$ is the solution}.

\noindent We would like to pass to the limit (\ref{4.2}) and obtain
\begin{equation}
\label{u}
u(x)=\int_D G_D^k(x,y)f(u(y))dy+b(x),\quad \text{for} \quad x\in
D,\quad u|_{\partial D}=0.
\end{equation}
By taking a subsequence of $u^{(n)}$ (still denoted by $u^{(n)}$), we can assume that $u^{(n)}$ converges to $u$, a.e., as $n\to \infty$.
Then, since $f_1$ is continuous $f_1(u^{(n)})\to f_1(u)$, a.e., as $n\to\infty$. Consequently, Schwarz inequality with respect to the measure
on $D$ with density $\vert G_D^k(x,y)\vert dy$ implies
\begin{align*}
&\left\Vert \int_D G_D^k(\cdot,y)\left[f_1(u^{(n)}(y))-f_1(u(y))\right] dy\right\Vert_{L^2(D)}^2\\
&\quad \le C \int_D dx \int_D dy  \left\vert G_D^k(x,y)\right\vert\left[f_1(u^{(n)}(y))-f_1(u(y))\right]^2 dy.
\end{align*}
By bounded convergence, this last expression goes to zero as $n\to\infty$.

Similarly,
\begin{align*}
&\left\Vert \int_D G_D^k(\cdot,y)\left[f_2(u^{(n)}(y))-f_2(u(y))\right] dy\right\Vert_{L^2(D)}^2\\
&\quad \le C \int_D dx \int_D dy  \left\vert G_D^k(x,y)\right\vert \vert u^{(n)}(y)-u(y)\vert^2\\
&\quad \le C \sup_{x\in D}\Vert G(x,\cdot)\Vert_{L^1(D)}\Vert u^{(n)}-u\Vert_{L^2(D)}^2,
\end{align*}
which clearly tends to zero as $n\to\infty$.

Thus, $u$ satisfies (\ref{u}) and
 $u\in \mathcal{S}$.
\hfill $\blacksquare $
\bigskip

The last part of the section is devoted to a further analysis of the solution $u$. In the next Lemma
we set $\tilde{\mathcal{K}}:= \sup_{x\in D} \int_D\left\vert G_D^k(x,y)\right\vert dy$.

\begin{lem}
\label{l4.3} Assume the  assumptions of Theorem $\ref{t4.1}$ and in addition, $L\tilde{\mathcal{K}}<1$.
Then, for any $p\in[1,\infty[$,
\begin{equation*}
\left\|u\right\|_{L^p(\Omega;L^{\infty}(D))}\leq C.
\end{equation*}
\end{lem}

\noindent{\it Proof}:  The  expression $b$ given in (\ref{defb})
belongs to $L^p(\Omega;L^{\infty}(D))$. Indeed, the function
$x\to\int_D G_D^k(x,y)g(y)dy$, defined on $D$, is continuous and
deterministic and therefore belongs to the space
$L^p(\Omega;L^{\infty}(D))$. As for $\int_D G_D^k(x,y)dB^H(y)$,
this property has been proved in Theorem \ref{t2.1}.

By definition, $u(x)= \int_D G_D^k(x,y)f(u(y))dy+b(x)$. Then, by the properties of $f$ we have
\begin{align*}
\left\|u\right\|_{L^p(\Omega;L^{\infty}(D))}&\le C\left\{1+L\left\|u\right\|_{L^p(\Omega;L^{\infty}(D))}\right\}
\sup_{x\in D} \int_D\left\vert G_D^k(x,y)\right\vert dy + \left\Vert b\right\Vert_{L^p(\Omega;L^{\infty}(D))}.
\end{align*}
Taking into account the restriction on the constant $L$, this yields the result
 \hfill \qed
\bigskip

We finally state a result on the regularity of the sample paths of the solution.
\begin{thm}
\label{t4.2} With the same hypotheses as in Theorem $\ref{t4.1}$, for any $x, z \in D$
and for any $p\in[1,\infty[$,
the solution $u$ to $(\ref{1.1})$ satisfies
\begin{equation}
\label{4.5} \mathbb{E}(|u(x)-u(z)|)^p\leq C|x-z|^{p\lambda},
\end{equation}
with $\lambda$ defined in Theorem \ref{t2.1} (see also Remark \ref{r3.1}).

Consequently, a.s. the sample paths are $\gamma$-H\"older
continuous with $\gamma\in]0,\lambda[$.
\end{thm}

{\it Proof}: We write $\left\Vert
u(x)-u(z)\right\Vert_{L^p(\Omega)}\leq \sum_{i=1}^3 I_{i}^k(x,z)$,
with
\begin{align*}
I_1^k(x,z)&=\left\Vert\int_D\left(G_D^k(x,y)-G_D^k(z,y)\right)f(u(y))dy\right\Vert_{L^p(\Omega)},\\
I_2^k(x,z)&=\left\Vert\int_D\left(G_D^k(x,y)-G_D^k(z,y)\right)g(y)dy\right\Vert_{L^p(\Omega)},\\
I_3^k(x,z)&=\left\Vert\int_D\left(G_D^k(x,y)-G_D^k(z,y)\right)dB^H(y)\right\Vert_{L^p(\Omega)}.
\end{align*}
By H\"older's inequality and $(\ref{holder})$ we have
\begin{align*}
I_1^k(x,z)&\leq
\left\|G_D^k(x,\cdot)-G_D^k(z,\cdot)\right\|_{L^{\frac{1}{H_1},\cdots,\frac{1}{H_k}}(D)}\\
&\quad\times\left(M+\left|f_2(0)\right|+L\left(\mathbb{E}\|u\|^p_{L^{\frac{1}{1-H_1},\cdots,\frac{1}{1-H_k}}(D)}\right)^\frac{1}{p}\right)\\
&\leq
C\left(M+\left|f_2(0)\right|+L\left(\mathbb{E}\|u\|^p_{L^{\frac{1}{1-H_1},\cdots,\frac{1}{1-H_k}}(D)}\right)^\frac{1}{p}\right)|x-z|^\lambda.
\end{align*}
The factor
$\mathbb{E}\|u\|^p_{L^{\frac{1}{1-H_1},\cdots,\frac{1}{1-H_k}}(D)}$
is finite. Indeed this follows from (\ref{2.17}) and Lemma
\ref{l4.3}. Consequently,

\begin{align*}
I_1^k(x,z)\leq C|x-z|^\lambda.
\end{align*}

In a similar but easier way, we obtain a similar bound for $I_2^k(x,z)$. As for $I_3^k(x,z)$, the bound is obtained by first applying the hypercontractivity
inequality and then (\ref{3.7}).

This ends the proof of (\ref{4.5}). The statement about the
regularity of the sample paths follows from Kolmogorov's
criterion. \hfill\qed


\section{Lattice approximations in $L^2$-spatial norm}
\label{s5} This section is devoted to give finite difference
approximation sequences for the SPDE (\ref{1.1}) on the domain
$D=]0,1[^k$, obtained by discretizing  the Laplacian operator.
To simplify the notation, we shall omit the superscript $k$
when referring to $G_D^k$, and denote by $C$ any positive, finite constant  not depending on $n$, although
its value may change from one line to the next.

The analysis of the rate of convergence is done by means of a general result given in the
next theorem.
\begin{thm}
\label{t5.frame}
Assume the hypotheses of Lemma \ref{l4.3}. We consider a sequence of functions $\{g_n\}_{n\ge 1}$
defined on $D$ such that
\begin{equation*}
\lim_{n\to\infty}\|g-g_n\|_{L^{\frac{1}{1-H_1},\ldots,\frac{1}{1-H_k}}(D)}=0,
\end{equation*}
along with another sequence of functions
 $\left\{\tilde G_{D,n}\right\}_{n\ge 1}$, defined on $D\times D$ satisfying
\begin{equation}
\label{greenn}
\int_D \left\Vert \tilde G_{D,n}(x,.)-G_D(x,.)\right\Vert_{L^{\frac{1}{H_1},\cdots,\frac{1}{H_k}}(D)}^2 dx \le C n^{-\gamma},
\end{equation}
for some $\gamma>0$.

Let  $\{\tilde u_n (x), x\in D\}_{n\ge 1}$ be a sequence of random fields satisfying
\begin{align}
\label{5.1.1}
\tilde u_n(x)=\int_D \tilde G_{D,n}(x,y)f(\tilde u_n(y))dy+\int_D \tilde G_{D,n}(x,y)g_n(y)dy+\int_D \tilde G_{D,n}(x,y)dB^H(y).
\end{align}
We suppose that for some $p\in[1,\infty[$ and $q\in\left[\frac{1}{1-\max(H_1,\ldots,H_k)},\infty\right]$,
\begin{equation}
\label{5.1.2}
\sup_{n\ge 1}\left\Vert \tilde u_n\right\Vert_{L^p(\Omega; L^q(D))}\le C.
\end{equation}
Then,
\begin{equation}
\label{5.1.3}
\left\Vert u-\tilde u_n\right\Vert_{L^p(\Omega; L^2(D))}\le C\left( n^{-\frac{\gamma}{4}}+\Vert g-g_n\Vert^{\frac{1}{2}}_{L^{\frac{1}{1-H_1},\ldots,\frac{1}{1-H_k}}(D)}\right),
\end{equation}
where $u$ is the solution of (\ref{1.3}).
\end{thm}
\noindent{\it Proof:}
We shall follow the scheme of the proof of
Theorem 2.4 of \cite{GM06} (see also Theorem 4 in \cite{MSS06}).
 By defining
\begin{align*}
T(x)&=\int_D[G_D(x,y)-\tilde G_{D,n}(x,y)]f(\tilde u_n(y))dy\notag\\
& +\int_DG_D(x,y)[g(y)-g_n(y)]dy
+\int_D[G_D(x,y)-\tilde G_{D,n}(x,y)]g_n(y)dy\notag\\
& +\int_D[G_D(x,y)-\tilde G_{D,n}(x,y)]dB^H(y),
\end{align*}
we have
$u(x)-\tilde u_n(x)=\int_D
G_D(x,y)[f(u(y))-f(\tilde u_n(y))]dy +T(x)$.

Using (\ref{2.0}) and (\ref{2.6.1}), as in \cite{MSS06} we obtain
\begin{align}
\label{5.10}
&(a-L)\|u-\tilde u_n\|_{L^2(D)}^2\nonumber\\
&\leq 2a\int_D(u(x)-\tilde u_n(x))T(x)dx
+\int_D\left(f(u(x))-f(\tilde u_n(x))\right)T(x) dx.
\end{align}

Let $q\in\left[\frac{1}{1-\max(H_1,\ldots,H_k)},\infty\right]$ and let $\tilde q\in[1,2]$ be its conjugate.
By applying H\"older's inequality and by virtue of the assumptions
on $f$, the right-hand side (\ref{5.10}) is bounded by
\begin{equation}
\label{5.11}
(2a+L)\|u-\tilde u_n\|_{L^q(D)}\|T\|_{L^{\tilde q}(D)}+ 2M\|T\|_{L^{\tilde q}(D)}.
\end{equation}
Lemma \ref{l4.3} yield
\begin{equation}
\label{5.12}
\sup_{n\ge1} \left\Vert u-\tilde u_n\right\Vert_{L^p(\Omega;L^q(D))}\le C,
\end{equation}
for the values of $p$ and $q$ such that (\ref{5.1.2}) holds true.

From (\ref{5.10})--(\ref{5.12}) and applying Schwarz's inequality 
we obtain
\begin{equation}
\left\Vert u-\tilde u_n\right\Vert_{L^p(\Omega;L^2(D))}\le C \left(\left\Vert T\right\Vert_{L^p(\Omega;L^{\tilde q}(D))}\right)^{\frac{1}{2}}.
\end{equation}

H\"older's inequality for $L^p$ spaces with mixed norm yields
\begin{align*}
|T(x)|&\le \|G_D(x,\cdot)-\tilde G_{D,n}(x,\cdot)\|_{L^{\frac{1}{H_1}, \cdots, \frac{1}{H_k}}(D)} \|f(\tilde u_n)\|_{L^{\frac{1}{1-H_1}, \cdots,\frac{1}{1-H_k}}(D)}\\
&\quad+\|G_D(x,\cdot)\|_{L^{\frac{1}{H_1},\cdots,\frac{1}{H_k}}(D)}\|g-g_n\|_{L^{\frac{1}{1-H_1},\cdots,\frac{1}{1-H_k}}(D)}\\
&\quad+\|G_D(x,\cdot)-\tilde G_{D,n}(x,\cdot)\|_{L^{\frac{1}{H_1},\cdots,\frac{1}{H_k}}(D)} \|g_n\|_{L^{\frac{1}{1-H_1},\cdots,\frac{1}{1-H_k}}(D)}\\
&\quad+\Big|\int_D[G_D(x,y)-\tilde G_{D,n}(x,y)]dB^H(y)\Big|.
\end{align*}

By the assumptions on the function $f$, we  have
\begin{equation*}
\|f(\tilde u_n)\|_{L^{\frac{1}{1-H_1},\cdots,\frac{1}{1-H_k}}(D)}
\le M+|f_2(0)|+L \|\tilde u_n\|_{L^{\frac{1}{1-H_1}, \cdots,\frac{1}{1-H_k}}(D)}.
\end{equation*}
Thus, by (\ref{2.17}) and  (\ref{5.1.2})
\begin{equation*}
\sup_{n\ge 1}\Vert f(\tilde u_n)\Vert_{L^p(\Omega;L^{\frac{1}{1-H_1},\cdots,\frac{1}{1-H_k}}(D))}
\le C.
\end{equation*}

Consequently, since $\tilde q \le 2$,
\begin{align}
\label{5.13}
\|T\|_{L^p(\Omega;L^{\tilde q}(D))}&\leq C\|T\|_{L^p(\Omega;L^2(D))}\notag\\
&\leq C\left\{n^{-\frac{\gamma}{2}}\Big(1+\|g\|_{L^{\frac{1}{1-H_1},\ldots,\frac{1}{1-H_k}}(D)}\Big)
+\Vert g-g_n\Vert_{L^{\frac{1}{1-H_1},\ldots,\frac{1}{1-H_k}}(D)}\right.\notag\\
&\left.+\left\Vert \int_D[G_D(\cdot,y)-\tilde G_{D,n}(\cdot,y)] dB^H(y)\right\Vert_{L^p(\Omega;L^2(D))}\right\},
\end{align}
where we have applied  (\ref{3.3.1}) and $\gamma$ is given in (\ref{greenn}).

Let us now give an upper bound for the last term in (\ref{5.13}). Assume first $p\in[1,2]$.
By applying H\"older's inequality with $\widetilde{p}=\frac{2}{p}\geq 1$
to the expectation operator and then Fubini's theorem, we obtain
\begin{align*}
&\left\Vert \int_D[G_D(\cdot,y)-\tilde G_{D,n}(\cdot,y)] dB^H(y)\right\Vert_{L^p(\Omega;L^2(D))}\\
&\quad=\Big(\mathbb{E}\Big\|\int_D[G_D(x,y)-\tilde G_{D,n}(x,y)]dB^H(y)\Big\|_{L^{2}(D)}^p\Big)^\frac{1}{p}\\
&\quad \leq \Big(\mathbb{E}\Big\|\int_D[G_D(x,y)-\tilde G_{D,n}(x,y)]dB^H(y)\Big\|_{L^{2}(D)}^{2}\Big)^\frac{1}{2}\\
&\quad\leq \Big(\int_D\Big(\mathbb{E}\Big|\int_D[G_D(x,y)-\tilde G_{D,n}(x,y)]dB^H(y)\Big|^2\Big)dx\Big)^\frac{1}{2}\\
&\quad \le \left( \int_D \left\Vert G_D(x,\cdot)-\tilde G_{D,n}(x,\cdot)\right\Vert_{L^{\frac{1}{H_1},\cdots,\frac{1}{H_k}}(D)}^2 dx\right)^{\frac{1}{2}}\\
&\quad\le C n^{-\frac{\gamma}{2}},
\end{align*}
where in the last inequality we have used (\ref{greenn}).

Next, we consider the case $p>2$.  We apply first Minkowski's inequality with
respect to the probability measure and the Lebesgue measure, then the
hypercontractivity inequality (see for instance \cite{LT91}), to obtain
\begin{align*}
&\left\Vert \int_D[G_D(\cdot,y)-\tilde G_{D,n}(\cdot,y)] dB^H(y)\right\Vert_{L^p(\Omega;L^2(D))}\\
&\quad =\Big(\mathbb{E}\Big(\int_D\Big|\int_D [G_D(x,y)-\tilde G_{D,n}(x,y)] dB^H(y)\Big|^{2}dx\Big)^\frac{p}{2}\Big)^\frac{1}{p}\\
&\quad=\Big\|\int_D\Big|\int_D [G_D(x,y)-\tilde G_{D,n}(x,y)]dB^H(y)\Big|^{2}dx\Big\|_{L^{\frac{p}{2}}(\Omega)}^{\frac{1}{2}}\\
&\quad\leq \Big(\int_D\Big\|\Big(\int_D [G_D(x,y)-\tilde G_{D,n}(x,y)]dB^H(y)\Big)^{2}\Big\|_{L^{\frac{p}{2}}(\Omega)}dx\Big)^\frac{1}{2}\\
&\quad=\Big(\int_D\Big(\mathbb{E}\Big|\int_D [G_D(x,y)-\tilde G_{D,n}(x,y)]dB^H(y)\Big|^p\Big)^\frac{2}{p}dx\Big)^\frac{1}{2}\\
&\quad\leq (p-1)^\frac{1}{2}\Big(\int_D\Big(\mathbb{E}\Big|\int_D [G_D(x,y)-\tilde G_{D,n}(x,y)]dB^H(y)\Big|^2\Big)dx\Big)^\frac{1}{2}\\
&\quad \le C(p)\left( \int_D \left\Vert G_D(x,\cdot)-\tilde G_{D,n}(x,\cdot)\right\Vert_{L^{\frac{1}{H_1},\cdots,\frac{1}{H_k}}(D)}^2 dx\right)^{\frac{1}{2}}\\
&\quad\le C(p) n^{-\frac{\gamma}{2}}.
\end{align*}

Plugging these estimates in (\ref{5.13}), we finish the proof of the
theorem. \hfill \qed
\bigskip

Let $I^k$ and $I_n^k$ be the sets of indices $\{1,2,\ldots \}^k$ and
$\{1,2,\ldots,n-1 \}^k$, respectively.
On the space  $X=\{ u:\; u=\{u_i\}_{i\in I^k_n}\}=\mathbb{R}^{(n-1)^k}$ endowed with the
Hilbert-Schmidt norm, we define the second order difference
operator $A:X\to X$,
\begin{equation*}
(Au)_i=\sum_{j=1}^k n^2\left(u_{i-e_j}-2u_i+u_{i+e_j}\right),\ i\in I_n^k,
\end{equation*}
where $\{e_j\}_{j=1}^k$ is the canonical basis of $\mathbb{R}^k$.

Consider the orthogonal complete system in $L^2(D)$ provided by the
functions
\begin{equation*}
v_{\beta}(x)=\sin(\beta_1\pi x_1)\cdots\sin(\beta_k\pi x_k),\
\beta\in I^k,\  k\geq 1.
\end{equation*}
The set of vectors
$\Big\{\Big(\frac{2}{n}\Big)^{k/2}U_\beta,\; \beta\in
I^k_n\Big\},\; (U_\beta)_i=v_{\beta}\Big(\frac{i}{n}\Big), \; i\in
I^k_n$,
is an orthonormal system in $X$ of eigenvectors of $A$, with
eigenvalues
$\lambda_\beta=-\pi^2(\beta_1^2c_{\beta_1}+\dots+\beta_k^2c_{\beta_k})$,
where
$c_l=\Big(\frac{l\pi}{2n}\Big)^{-2}\sin^2\Big(\frac{l\pi}{2n}\Big)$.
Notice that $\frac{4}{\pi^2}\leq c_l\leq 1$ for every $1\leq l\leq
n-1$.
\medskip

We will consider approximation schemes of (\ref{1.3}) based on the
grid of  $\bar D$ given by
\begin{equation*}
\mathcal{G}=\Big\{\frac{j}{n}=\Big(\frac{j_1}{n},\ldots,
\frac{j_k}{n}\Big):\ j_l=0,1,\ldots,n,\; l=1,\ldots,k\Big\}\subset
\bar D.
\end{equation*}
For any point $\frac{j}{n}\in\mathcal{G}$, we set
$D_j=[\frac{j_1}{n},\frac{j_1+1}{n}[\times\cdots\times[\frac{j_k}{n}, \frac{j_k+1}{n}[$, and for each $x\in D_j$
 we define $\kappa_n(x)=\frac{j}{n}$.

We begin by giving a first type of discrete approximations of $u$ on points of
$\mathcal{G}$, denoted by $u_n$. If $\frac{j}{n}\in
\mathcal{G}\cap\partial D$, we set
$u_n(\frac{j}{n})=0$ (boundary conditions), while for $\frac{j}{n}$ with $j\in I_n^k$, we define
$u_n(\frac{j}{n})$ to be the solution of the system
\begin{equation}
\label{k1}
A u_n=f(u_n)+g_n+n^k\mathbf{B}^H,
\end{equation}
where $\mathbf{B}^H$ is the vector $\{B^H(D_i)=
\int_{\mathbb{R}^k}\Ind_{D_i}(y)dB^H(y),\; i\in I_n^k\}$ and $(g_n)_{n\ge 1}$ a sequence of step functions,
$g_n(x)=g_n(\kappa_n(x))$, $n\ge 1$. Then, for any $x\in D$ we
define $u_n(x)=u_n(\kappa_n(x))$. From \cite{GM06} we know that
$\{u_n(x), x\in D\}$ satisfies the evolution equation
\begin{equation}
\label{k2}
u_n(x)=\int_D G_{D,n}(x,y) f(u_n(y)) dy + \int_D G_{D,n}(x,y) g_n(y) +\int_D G_{D,n}(x,y)dB^H(y),
\end{equation}
where
\begin{equation}
\label{k3}
G_{D,n}(x,y)=\sum_{\beta\in I_n^k} \frac{2^k}{\lambda_\beta} v_\beta(\kappa_n(x))v_\beta(\kappa_n(y)).
\end{equation}
\medskip

In dimension $k=1,2,3$, we shall consider $\{u_n(x), x\in D\}$, $n\ge 1$, as sequence of approximations of
the process $\{u(x), x\in D\}$. We notice that in this case the kernel $G_{D,n}(x,.)$ is related with the truncation of
the Fourier expansion of $G_D(x,.)$. For $k\ge 4$ we shall follow the more sophisticated approach of \cite{MSS06},
which considers a smoothed version of  $G_D(x,.)$. We remark that for such dimensions $G_D(x,.)$ is not square integrable.

For low dimensions, we have the following.
\begin{thm}
\label{t5.low}  Assume the assumptions of Lemma \ref{l4.3}. Let
$k\in\{1,2,3\}$ and $\{u_n(x), x\in D\}$, $n\ge 1$, be defined in
(\ref{k2}), (\ref{k3}). Then, for any $p\in[1,\infty[$,
\begin{equation}
\label{conv1}
 \|u-u_n\|_{L^p(\Omega; L^2(D))}\leq
C \left(n^{-\nu} + \|g-g_n\|^{\frac{1}{2}}_{L^{\frac{1}{1-H_1},\cdots,\frac{1}{1-H_k}}(D)}\right),
\end{equation}
with $\nu\in]0,\frac{1}{2}]$,  $\nu\in]0,\frac{1}{2}[$, $\nu\in]0,\frac{1}{4}]$, for
$k=1$, $k=2$ and $k=3$, respectively.
\end{thm}

\noindent{\it Proof:} It follows from Theorem \ref{t5.frame}. Indeed, first we give a slight improvement
of Lemma 3.4 in \cite{GM06} which yields the validity of
condition (\ref{greenn}) for
$\tilde G_{D,n}:=G_{D,n}$ with $\gamma:=4\nu$ and the
values of $\nu$ of the statement.

In fact, for the expression termed $A$ in \cite{GM06}, page 223, we have
\begin{equation*}
A\le C \int_n^\infty r^{k-1-4} dr = Cn^{k-4},
\end{equation*}
while for $k=3$, we can proceed with the term called $B$ as follows. Let $\rho\in]1,2[$, then
\begin{align*}
B&\le C\sum_{\alpha\in I_n^k}\frac{1}{|\alpha|^2n^2} = \frac{C}{n^\rho}\sum_{\alpha\in I_n^k}\frac{1}{|\alpha|^{4-\rho}}\\
&\le C n^{-1}.
\end{align*}

The process $\tilde u_n:=u_n$ satisfies (\ref{5.1.2}). This can be easily checked by applying Lemma 3.3 in \cite{GM06}
and similar arguments as those in the proof of Lemma \ref{l5.4} below  (we leave the details to the reader).
\hfill\qed
\medskip

We next deal with higher dimensions.
 The Fourier analysis techniques we shall use in the proofs require
the identification of functions $f:[-1,1[^k\longrightarrow{\mathbb{R}}$ with functions
$F:\mathbb{T}^k\longrightarrow{\mathbb{R}}$ defined on the $k$-th dimensional torus through the exponential mapping
$F(e^{i\pi x}):=F(e^{i\pi x_1},\ldots,e^{i\pi x_k})$,
which carries Lebesgue measure into the Haar measure, that is
$\int_{]-1,1[^k}f(x)dx=\int_{\mathbb{T}^k}F(e^{i\pi x})dx$.
When dealing with the function $y\to G_D(x,y)$, we will consider
its odd extension,
that is, for any $x, y\in]0,1[^k$, $y=(y_1,\ldots,y_k)$, we define
$G_D(x,(y_1,\dots,-y_i,\dots,y_k))=-G_D(x,(y_1,\dots,y_i,\dots,y_k))$.
We still note $G_D(x,\cdot)$ the extension. Let
$\mathbb{G}_D^{x}(e^{i\pi y})=G_D(x,y)$
be its identification with a function defined on $\mathbb{T}^k$. Observe that
$\mathbb{G}_D^{x}$ satisfies
\begin{equation*}
\|\mathbb{G}_D^{x}(e^{i\pi\cdot})\|_{L^{p_1,\ldots,p_k}(\mathbb{T}^k)}
=\|G_D(x,\cdot)\|_{L^{p_1,\ldots,p_k}(]-1,1[^k)}=2^{\sum_{j=1}^k\frac{1}{p_j}}\|G_D(x,\cdot)\|_{L^{p_1,\ldots,p_k}(D)},
\end{equation*}
for any $p_1,\ldots,p_k$ such that the last norm is finite.

Let $\psi(x)\in{\mathcal C}^\infty_c(]-1,1[)$ be an even function,
$0\leq \psi\leq 1$, $\int_{-1}^1\psi=1$. Set
$\Psi(x)=\prod_{j=1}^k\psi(x_j)$. Clearly,
$\Psi(x)\in{\mathcal C}^\infty_c(]-1,1[^k)$ and it is an even
function in each variable $x_j$. Define
\begin{equation*}
\Phi(e^{i\pi x})=\prod_{j=1}^k\phi(e^{i\pi x_j}) :=
\prod_{j=1}^k\psi(x_j)=\Psi(x).
\end{equation*}
The functions
$\Phi_\varepsilon(e^{i\pi x})=\frac1{\varepsilon^k}\Psi\big(\frac{x}{\varepsilon}\big):=\Psi_\varepsilon(x)$,
$\varepsilon>0$, provide an approximation of the identity in
$\mathbb{T}^k$.

We shall denote by $\hat\Psi$ the Fourier transform of $\Psi$, which is a rapidly decreasing function, therefore for any  $\theta\in[0,\infty[$
there is a constant $C(\theta)$ such that $\sup_{\xi}|\xi|^\theta|\hat\Psi(\xi)|\le C(\theta)$.
\medskip

Let us now introduce a second kind of approximations of $u$. For this we start by writing
 $A=U^t\Lambda U$,  with
$U$  the $(n-1)^k$ matrix whose rows are the vectors
$U_{\beta_j}$, (here $\beta_j$, $j=1,\cdots,(n-1)^k$ is the
lexicographic enumeration of $I_n^k$) and  $\Lambda$ the square diagonal matrix with entries
$\Lambda_{j,j}=\lambda_{\beta_j}$.

The smoothed version of $A$ is defined as follows.
Fix $\varepsilon>0$ and  define $\Lambda^{\varepsilon}$ as the square
diagonal matrix in dimension $(n-1)^k$ with diagonal entries
$$\lambda^\varepsilon_{\beta_j}=\frac{\lambda_{\beta_j}}{\hat\Psi(\varepsilon\beta_j)}.$$
In connection with $\Lambda^\varepsilon$ we define a sequence $(u_n^\varepsilon, n\ge 1)$ of functions in the following way.
If $\frac{j}{n}\in
\mathcal{G}\cap\partial D$, set
$u_n^\varepsilon(\frac{j}{n})=0$ (boundary conditions). For
$\frac{j}{n}$, with $j\in I_n^k$, define
$u_n^\varepsilon(\frac{j}{n})$ to be the solution of the system
\begin{equation}
\label{5.4}
(U^t \Lambda^\varepsilon U)
u_n^\varepsilon=f(u_n^\varepsilon)+g_n+n^k\mathbf{B}^H.
\end{equation}
Finally, for any $x\in D$ we
define $u_n^\varepsilon(x)=u_n^\varepsilon(\kappa_n(x))$.

We shall prove later that an appropriate sequence
$u_n:=u_n^{\varepsilon(n)}$ of such functions converges to
the solution of (\ref{1.1}) in the space $L^p(\Omega;L^2(D))$,
for any $p\ge 1$, with a rate of convergence which depends on the dimension
$k$, on the driving noise and on the rate of convergence of $g_n$ to $g$.
\medskip

The following result is proved with the same arguments as in Proposition 1 of
\cite{MSS06}.
\begin{prop}
\label{p5.1} With the same hypotheses as in Theorem \ref{t4.1}
and assuming that the Lipschitz constant satisfies $L<4k$,
Equation (\ref{5.4}) possesses a unique solution.
Moreover, this solution satisfies the {\it mild} equation
\begin{align}
\label{5.5}
u_n^\varepsilon(x)&=\int_DG_{D,n}^\varepsilon(x,y)f(u_n^\varepsilon(y))dy
+\int_DG_{D,n}^\varepsilon(x,y)g_n(y)dy+\int_DG_{D,n}^\varepsilon(x,y)dB^H(y),
\end{align}
where
\begin{equation}
\label{5.6}
G_{D,n}^\varepsilon(x,y)=\sum_{\beta\in
I^k_n}\frac{\hat\Psi(\varepsilon\beta)2^k}
{\lambda_{\beta}}v_{\beta}(\kappa_n(x))v_{\beta}(\kappa_n(y)).
\end{equation}
\end{prop}

Both (\ref{k3}) and (\ref{5.6})  correspond to  discretized  Fourier
series expansions; in (\ref{5.6}), the Fourier coefficients are smoothed by the factor $\hat\Psi(\varepsilon\beta)$.
\medskip

Our next aim is to apply Theorem \ref{t5.frame}  to $\tilde u_n:=u_n^{\varepsilon(n)}$ defined in (\ref{5.5}), for values of
$\varepsilon$ that depend on $n$, and dimensions $k\ge 4$.
The next statements provide the ingredients for checking  condition (\ref{greenn}) for $\tilde G_{D,n}:= G_{D,n}^{\varepsilon(n)}$.
\medskip

Set $G_D^\varepsilon(x,y)=\mathbb{G}_D^{x,\varepsilon}(e^{i\pi y})$, where
$\mathbb{G}_D^{x,\varepsilon}(e^{i\pi y})=
\int_{\mathbb{T}^k}\mathbb{G}_D^{x}(e^{i\pi
(y-u)})\Phi_\varepsilon(e^{i\pi u})du$.
The function $G_D^\varepsilon(x,y)$ is a smoothing of $G_D(x,\cdot)$.

\begin{lem}
\label{l5.2}
For any $\varepsilon>0$, we have
\begin{equation}
\label{5.1}
G_D^\varepsilon(x,y)=\sum_{\beta\in
I^k}\frac{-\hat{\Psi}(\varepsilon\beta)2^k}{\pi^2|\beta|^2}
v_{\beta}(x)v_{\beta}(y),
\end{equation}
in $L^2(D\times D)$ and a.e.
In addition,
\begin{equation}
\label{5.2} \|G_D^{\varepsilon}(x,\cdot)\|_{L^2(D)}^2
=\frac{2^k}{\pi^4}\sum_{\beta \in
I^k}\frac{\hat{\Psi}^2(\varepsilon\beta)}{|\beta|^4}
v_{\beta}^2(x),
\end{equation}
and the series converges uniformly in $x\in D$ and $\varepsilon\in]0,\varepsilon_0]$.
\end{lem}

\noindent{\it Proof}: The first part of the assertion is Lemma 8 of \cite{MSS06}. For the sake of completeness and further use, we
give some details of its proof. Fix $p\in\left[1,\frac{k}{k-2}\right[$. Young's inequality for convolution  (\cite{AF03}, page 34, Corollary 2.25) yields

\begin{align*}
\|G_D^{\varepsilon}(x,\cdot)\|_{L^2(D)} &=
2^{-\frac{k}{2}}\|\mathbb{G}_D^{x,\varepsilon}(e^{i\pi \cdot})\|_{L^2(\mathbb{T}^k)}\\
&\leq
2^{-\frac{k}{2}}\|\mathbb{G}_D^x\|_{L^{p}(\mathbb{T}^k)}\|\Phi_\varepsilon\|_{L^{r}(\mathbb{T}^k)}\\
&= C\|G_D(x,\cdot)\|_{L^{p}(D)}\|\Psi_\varepsilon\|_{L^{r}(]-1,1[^k)},
\end{align*}
for $\frac{1}{2}=\frac{1}{r}+\frac{1}{p}-1$.

For any $r\ge 1$,
$\sup_\varepsilon\Vert\Psi_\varepsilon\Vert_{L^{r}(]-1,1[^k)}\le C$. Thus, using (\ref{3.1}) we obtain
\begin{equation}
\label{5.3} \sup_{x\in D}\sup_{\varepsilon\in]0,\varepsilon_0]}\|G_D^{\varepsilon}(x,\cdot)\|_{L^2(D)}\le C.
\end{equation}
Hence $G_D^\varepsilon\in L^2(D\times D)$ and the formula (\ref{5.1}) follows from the computation of the Fourier coefficients carried out in \cite{MSS06}, Lemma 8.

The orthogonal complete
system $(v_\beta,\beta\in I^k)$ satisfies $\|v_{\beta}\|_{L^2(D)}=2^{-k/2}$. Thus (\ref{5.2}) follows easily from (\ref{5.1}). Finally,
 $(\ref{5.3})$ implies the uniform convergence of the series
 in (\ref{5.2}). \hfill \qed
\bigskip

The next result provides an estimate of the discrepancy between $G_D$ and
$G_D^\varepsilon$.

\begin{lem}
\label{l5.1} Assume  {\bf(${\rm \bf H}$)}. Then for every
$\varepsilon>0$,
\begin{equation}
\label{norm}
\sup_{x\in
D}\|G_D(x,\cdot)-G_D^\varepsilon(x,\cdot)\|_{L^{\frac{1}{H_1},\ldots,\frac{1}{H_k}}(D)}\leq
C\varepsilon^\lambda,
\end{equation}
with
\begin{equation*}
\left\{\renewcommand{\arraystretch}{1.8}
\begin{array}{ll}
\lambda=1,&\text{ for
}k=1,\\
\lambda\in\left]0,\left(2-k+\sum_{i=1}^k H_i\right)\wedge 1\right[, &\text{
for }k \ge 2.
\end{array}
\right.
\end{equation*}
\end{lem}

\noindent{\it Proof}: Since $G_D(x,\cdot)$, $G_D^\varepsilon(x,\cdot)$ are
odd in the $y_i$-variables, we have
\begin{align*}
\|G_D(x,\cdot)-G_D^\varepsilon(x,\cdot)\|_{L^{\frac{1}{H_1},\ldots,\frac{1}{H_k}}(D)}
&=
C\|G_D(x,\cdot)-G_D^\varepsilon(x,\cdot)\|_{L^{\frac{1}{H_1},\ldots,\frac{1}{H_k}}(]-1,1[^k)}\\
&=
C\|\mathbb{G}_D^x(e^{i\pi\cdot})-\mathbb{G}_D^{x,\varepsilon}(e^{i\pi\cdot})\|_{L^{\frac{1}{H_1},\ldots,\frac{1}{H_k}}(\mathbb{T}^k)},
\end{align*}
with $C=2^{-\sum_{i=1}^kH_i}$. Using that
$\int_{\mathbb{T}^k}\Phi_\varepsilon(e^{i\pi u})du=1$, we can
write
\begin{align*}
&\|\mathbb{G}_D^x(e^{i\pi\cdot})-\mathbb{G}_D^{x,\varepsilon}(e^{i\pi\cdot})\|_{L^{\frac{1}{H_1},\ldots,\frac{1}{H_k}}(\mathbb{T}^k)}\\
&\quad=\Big\|\int_{\mathbb{T}^k}(\mathbb{G}_D^x(e^{i\pi\cdot})-\mathbb{G}_D^x(e^{i\pi(\cdot-u)}))\Phi_\varepsilon(e^{i\pi u})du\Big\|_{L^{\frac{1}{H_1},\ldots,\frac{1}{H_k}}(\mathbb{T}^k)}\\
&\quad\le
\int_{\mathbb{T}^k}\|\mathbb{G}_D^x(e^{i\pi\cdot})-\mathbb{G}_D^x(e^{i\pi(\cdot-u)})\|_{L^{\frac{1}{H_1},\ldots,\frac{1}{H_k}}(\mathbb{T}^k)}\Phi_\varepsilon(e^{i\pi u})du\\
&\quad= 2^{\sum_{i=1}^kH_i}
\int_{\mathbb{T}^k}\|G_D(x,\cdot)-G_D(x,\cdot-u)\|_{L^{\frac{1}{H_1},\ldots,\frac{1}{H_k}}(D)}\Phi_\varepsilon(e^{i\pi u})du\\
&\quad\le C \int_{\mathbb{T}^k} |u|^\lambda \Phi_\varepsilon(e^{i\pi u})du\\
&\quad= C\int_{]-1,1[^k}|x|^{\lambda}\frac{1}{\varepsilon^k}\Psi\left(\frac{x}{\varepsilon}\right) dx,
\end{align*}
where  we have applied Minkowski's inequality with respect to the
finite measure on $\mathbb{T}^k$ defined by
$\Phi_\varepsilon(u)du$ and eventually (\ref{holder}). From this we infer (\ref{norm}) by observing that
the support of the function $\Psi$ is included in $]-1,1[^k$.

\hfill \qed
\medskip

As an additional auxiliary result, we need {\it a priori estimates} for the solution of (\ref{5.5}). An ingredient for this
is provided by the following Lemma.

\begin{lem}
\label{l5.3}
Fix $\varepsilon_0>0$. The smoothed, discretized Green function defined in (\ref{5.6}) satisfies
\begin{equation}
\label{5.7}
\sup_{n\ge 1}\sup_{x\in D, \varepsilon\in]0,\varepsilon_0]}\|G_{D,n}^\varepsilon(x,\cdot)\|_{L^2(D)}< +\infty.
\end{equation}
\end{lem}
{\it Proof}: The system $\{v_{\beta}(\kappa_n(y))\}$ is orthogonal  in $\mathbb{R}^{(n-1)^k}$, thus
in $L^2(D)$ as well.
Hence, using the lower bound $|\lambda_\beta|\geq 4|\beta|^2$ we have
\begin{equation*}
\|G_{D,n}^\varepsilon(x,\cdot)\|_{L^2(D)}^2 = \sum_{\beta\in
I^k_n}\frac{\hat{\Psi}^2(\varepsilon\beta)2^{k}}{\lambda_{\beta}^2}v_\beta^2(\kappa_n(x))
\leq \frac{2^{k}}{16} \sum_{\beta\in I^k_n}\frac{\hat{\Psi}^2(\varepsilon\beta)}{|\beta|^4}v_\beta^2(\kappa_n(x)).
\end{equation*}
From the last inequality, along with (\ref{5.2}), we obtain
\begin{align*}
\sup_{n\ge 1}\sup_{x\in D, \varepsilon\in]0,\varepsilon_0] }\|G_{D,n}^\varepsilon(x,\cdot)\|_{L^2(D)}^2
&\leq C\sup_{x\in D, \varepsilon\in]0,\varepsilon_0]}\Big\{\sum_{\beta\in I^k}
\frac{\hat{\Psi}^2(\varepsilon\beta)}{|\beta|^4}v_\beta^2(x)\Big\}\\
&\le C \sup_{x\in D, \varepsilon\in]0,\varepsilon_0]}\left\Vert G_D^\varepsilon(x,\cdot)\right\Vert_{L^2(D)}^2 <+\infty.
\end{align*}
\hfill \qed

Lemma 9 of \cite{MSS06} gives a more particular statement than the previous Lemma \ref{l5.3}. We have found an incorrect
argument in the proof of  the former that can be fixed using the proof of the later.
\medskip

Let $\bar q\in[1,2]$. By (\ref{2.17}), H\"older's inequality and Lemma \ref{l5.3} we have
\begin{equation}
\label{constant}
\sup_{n\ge 1}\sup_{x\in D,\varepsilon\in]0,\varepsilon_0]}
\left(\|G_{D,n}^\varepsilon(x,\cdot)\|_{L^{\frac{1}{H_1},\cdots,\frac{1}{H_k}}(D)}\vee
\|G_{D,n}^\varepsilon(x,\cdot)\|_{L^{\bar q}(D)}\right)\le \widehat{\mathcal{K}},
\end{equation}
for some positive, finite constant $\widehat{\mathcal{K}}$.
\medskip

We can now prove an {\it a priori} estimate for the solution of (\ref{5.5}).

\begin{lem}
\label{l5.4} Assume the same assumptions as in Proposition \ref{p5.1}
 with the Lipschitz constant  satisfying the restriction
 $L< \min(4k,\widehat{\mathcal{K}}^{-1})$, where $\widehat{\mathcal{K}}$ is given
in $(\ref{constant})$. Then, for any $p\in[1,\infty[$ and $q\in[2,\infty[$,
\begin{equation*}
\sup_{n\geq 1}\sup_{\varepsilon\in]0,\varepsilon_0]}
\left(\|u^{\varepsilon}_n\|_{L^p(\Omega;L^{q}(D))}\right)\leq C.
\end{equation*}
\end{lem}
{\it Proof}: Fix $q\in[2,\infty[$ and denote by  $\bar q\in]1,2]$ its  conjugate. H\"older's inequality and the properties on $f$ imply, for any $x\in D$,
\begin{align}
\label{fl5.4}
|u^\varepsilon_n(x)|&\leq\sup_{n\ge 1}\sup_{x\in D, \varepsilon\in]0,\varepsilon_0]}
\|G_{D,n}^\varepsilon(x,\cdot)\|_{L^{\bar q}(D)}
\Big(M+|f_2(0)|+L\|u^\varepsilon_n\|_{L^{q}(D)}\Big)\nonumber\\
&+\sup_{n\ge 1}\sup_ {x\in D, \varepsilon\in]0,\varepsilon_0]}
\|G_{D,n}^\varepsilon(x,\cdot)\|_{L^{\frac{1}{H_1},\cdots,\frac{1}{H_k}}(D)}\sup_{n\ge1}\Vert
g_n\Vert_{L^{\frac{1}{1-H_1},\ldots,\frac{1}{1-H_k}}(D)}\nonumber\\
&+ \Big|\int_D G_{D,n}^\varepsilon(x,y)dB^H(y)\Big|
\end{align}

Since $u_n^\varepsilon$ is a step function, its
$L^q$--norm is finite.
Moreover, arguing in a similar manner as we did in Theorem \ref{t5.frame} to get an upper bound on the last term of (\ref{5.13}), and applying (\ref{constant}) we obtain
\begin{equation*}
\Big\|\int_D G_{D,n}^\varepsilon(\cdot,y)
dB^H(y)\Big\|_{L^p(\Omega;L^q(D)}
\le C(p,q,H).
\end{equation*}
Taking the $L^p(\Omega;L^q(D)$--norm in  (\ref{fl5.4}) yields the conclusion.
\hfill \qed
\medskip

The next lemma gives an estimate of the discrepancy between
$G_D^\varepsilon$ and $G_{D,n}^\varepsilon$.
\smallskip

\begin{lem}
\label{l5.5} Assume $k\ge 4$. Let $\delta\in]0,2[$,  $\mu\in\left]0,\frac{2-\delta}{k-2}\right[$ and set
$\varepsilon(n)=n^{-\mu}$. There exists a positive, finite constant
$C(\delta,k)$, depending on $\delta$ and $k$ but not on $n$, such that
\begin{equation}
\label{5.8}
\|G_D^{\varepsilon(n)}-G_{D,n}^{\varepsilon(n)}\|_{L^2(D\times D)}\leq C(\delta,k) n^{-\frac{\delta}{2}}.
\end{equation}
\end{lem}

\noindent{\it Proof}: We follow the proofs of Lemma 3.4 in \cite{GM06} and Lemma 10 in
\cite{MSS06}.  By the definitions of the kernels
$G_{D,n}^\varepsilon$ and $G_D^\varepsilon$ given in (\ref{5.6})
and (\ref{5.1}) respectively, we have
\begin{equation*}
\|G_D^{\varepsilon(n)}-G_{D,n}^{\varepsilon(n)}\|_{L^2(D\times D)}^2
\le C\sum_{i=1}^4 A_i(x,y),
\end{equation*}
with
\begin{align*}
A_1&=\int_{D\times D}\Big|\sum_{\beta\in I^k\setminus I^k_n}\frac{-2^k\hat{\Psi}(\varepsilon\beta)}{\pi^2|\beta|^2} v_{\beta}(x)v_{\beta}(y)\Big|^2dx dy,\\
A_2&=\int_{D\times D}\Big| \sum_{\beta\in I^k_n}
\Big[\frac{-1}{\pi^2|\beta|^2}-\frac{1}{\lambda_\beta}\Big]
2^k\hat{\Psi}(\varepsilon\beta)
v_{\beta}(x)v_{\beta}(y)\Big|^2dx dy,\\
A_3&=\int_{D\times D}\Big| \sum_{\beta\in I^k_n}
\frac{2^k\hat{\Psi}(\varepsilon\beta)}{\lambda_\beta}
\left[v_{\beta}(x)-v_{\beta}(\kappa_n(x))\right]v_{\beta}(y)\Big|^2dx dy,\\
A_4&=\int_{D\times D}\Big| \sum_{\beta\in I^k_n}
\frac{2^k\hat{\Psi}(\varepsilon\beta)}{\lambda_\beta}
v_{\beta}(\kappa_n(x))\left[v_{\beta}(y)-v_{\beta}(\kappa_n(y))\right]\Big|^2dx
dy.
\end{align*}

In the sequel, we shall write $\varepsilon$ instead of
$\varepsilon(n)$ for simplicity, and we fix $\delta>0$.
Remember that $(v_\beta, \beta\in I^k)$ is a family of orthogonal functions in $L^2(D)$
with $\Vert v_\beta\Vert_{L^2(D)}=2^{-\frac{k}{2}}$.
Let $\theta>\frac{k-4}{2}$. Since $\hat{\Psi}$ is a rapidly decreasing function, we have
\begin{align*}
A_1&=\sum_{\beta\in I^k\setminus I^k_n}\frac{\hat{\Psi}^2(\varepsilon \beta)}{\pi^4|\beta |^4}
 \leq\frac{C(\theta)}{\varepsilon^{2\theta}}\sum_{\beta\in I^k\setminus I^k_n}\frac{1}{|\beta|^{4+2\theta}}\\
&\le {C(\theta)} \varepsilon^{-2\theta} n^{-4-2\theta+k} = {C(\theta)} n^{2\theta\mu-4-2\theta+k}.
\end{align*}
Fix $\theta:=\frac{\delta+k-4}{2-2\mu}$ in the last expression. We obtain
$A_1\le C(\theta) n^{-\delta}$.

For the analysis of the term $A_2$ we apply the estimate
$\Big|\frac{-1}{\pi^2|\beta|^2}-\frac{1}{\lambda_\beta}\Big|\le
\frac{C}{|\beta| n}$, valid for any $\beta\in I^k$. Taking $\theta<\frac{k-2}{2}$, we obtain
\begin{align*}
A_2& = \sum_{\beta\in I^k_n}\Big|\frac{-1}{\pi^2|\beta|^2}-\frac{1}{\lambda_\beta}\Big|^2\hat{\Psi}^2(\varepsilon\beta)
\leq  \frac{C}{n^2}\sum_{\beta\in I^k_n}\frac{\hat{\Psi}^2(\varepsilon\beta)}{|\beta|^2}\\
&\leq \frac{C(\theta)}{\varepsilon^{2\theta} n^2}
\sum_{\beta\in I^k_n}\frac{1}{|\beta|^{2+2\theta}}\le C(\theta) n^{2\mu\theta-4+k-2\theta}.
\end{align*}
Consider $\theta:=\frac{k-4+\delta}{2(1-\mu)}$. The last estimates yield
$A_2\le C(\theta) n^{-\delta}$.

For the study of the remaining terms, we use that for any $\beta\in I^k$,
$|v_{\beta}(x)-v_{\beta}(z)|\leq C|\beta|\,|x-z|$, and $|\lambda_\beta|\ge 4|\beta|^2$. This ensures
\begin{equation*}
\max\left(A_3, A_4\right)\le  \frac{C}{n^2}\sum_{\beta\in I^k_n}\frac{\hat{\Psi}^2(\varepsilon\beta)|\beta|^2}{\lambda_\beta^2}
\leq \frac{C(\theta)}{\varepsilon^{2\theta} n^2}
\sum_{\beta\in I^k_n}\frac{1}{|\beta|^{2+2\theta}}.
\end{equation*}
Hence, as for $A_2$, we obtain $\max\left(A_3, A_4\right)\le C(\theta) n^{-\delta}$.

The proof of the Lemma is now complete.
 \hfill \qed

As a consequence of Lemma \ref{l5.1}, (\ref{2.17}) and
Lemma \ref{l5.5} we obtain the following result.

\begin{cor}
\label{c5.1} With the same assumptions as in Lemmas \ref{l5.1} and \ref{l5.5},
there exists a positive constant $C$ not
depending on $n$, such that
\begin{equation}
\label{5.9}
\int_D\|G_D(x,\cdot)-G_{D,n}^{\varepsilon(n)}(x,\cdot)\|_{L^{\frac{1}{H_1},\cdots,\frac{1}{H_k}}(D)}^2dx\leq
C n^{-{\gamma}},
\end{equation}
with $\gamma=(2\mu \lambda) \wedge \delta$.
\end{cor}
\smallskip

Assume {\bf(${\rm \bf H}^{\ast}$)}. Then in the preceding
Corollary, $\lambda=1$ and $\gamma=2\mu\wedge \delta$. The biggest
value of $\gamma$ occurs when $2\mu=\delta$. Since
$\mu<\frac{2-\delta}{k-2}$ and the equation
$\frac{2(2-\delta)}{k-2}=\delta$ has the solution
$\delta=\frac{4}{k}$, we conclude that
$\gamma\in\left]0,\frac{4}{k}\right[$.

Suppose next that $k-1\geq\sum_{i=1}^k H_i>k-2$.  In this case
$$
2\mu\lambda\in\left]0,\frac{(4-2\delta)\left(2-k+\sum_{i=1}^kH_i\right)}{k-2}\right[.
$$
As before, the biggest upper bound of $\gamma$ is obtained by
solving the equation
$\frac{(4-2\delta)\big(2-k+\sum_{i=1}^kH_i\big)}{k-2}=\delta$,
which leads to the value $\delta=
\frac{4\big(2-k+\sum_{i=1}^kH_i\big)}{2-k+2\sum_{i=1}^kH_i}$.
Therefore
$\gamma\in\left]0,\frac{4\big(2-k+\sum_{i=1}^kH_i\big)}{2-k+2\sum_{i=1}^kH_i}\right[$.
\medskip

Consequently, under the assumption {\bf (H)}, we obtain
\begin{equation*}
\gamma\in\left]0,\frac{4\left(2-k+\left(\sum_{i=1}^kH_i\right)\wedge (k-1)\right)}{2-k+2\left(\left(\sum_{i=1}^kH_i\right)\wedge(k-1)\right)}\right[.
\end{equation*}


We have now the ingredients to  establish the convergence of the discretized scheme defined in (\ref{5.5}) to the solution of (\ref{1.1}) when
the parameters $n$ and $\varepsilon$ are related by the same constraints as in Lemma \ref{l5.5}. The constant 
$ \widehat{\mathcal{K}}$ in the next statement is  given in 
(\ref{constant}).

\begin{thm}
\label{t5.1} We assume the hypotheses of Lemmas \ref{l4.3} and in addition, 
$L<\min(4k, \widehat{\mathcal{K}}^{-1})$. 
 Let $k\ge 4$ and $\varepsilon(n)$ be as in Lemma
\ref{l5.5}. Then for any $p\in[1,\infty[$,
\begin{equation}
\label{conv2}
 \|u-u^{\varepsilon(n)}_n\|_{L^p(\Omega; L^2(D))}\leq
C \left(n^{-\nu}+ \Vert g-g_n\Vert^{\frac{1}{2}}_{L^{\frac{1}{1-H_1}, \ldots, \frac{1}{1-H_k}}(D)}\right),
\end{equation}
with
$\nu\in\left]0,\frac{2-k+\left(\sum_{i=1}^kH_i\right)\wedge (k-1)}{2-k+2\left(\left(\sum_{i=1}^kH_i\right)\wedge(k-1)\right)}\right[$.
\end{thm}

\noindent{\it Proof}: It follows from Theorem \ref{t5.frame}, 
Corollary \ref{c5.1} and the discussion that precedes the statement.
\hfill\qed
\medskip


Let $p_0:= \max_{i=1,\ldots,k}\left(\frac{1}{1-H_i}\right)$ and assume that
\begin{equation}
\label{rate}
\Vert g-g_n\Vert_{L^{p_0}(D)}\le C n^{-1},
\end{equation}
with a constant $C$ independent of $n$. Then the right hand side of (\ref{conv1}) and (\ref{conv2}) can be replaced by
$Cn^{-\nu}$. We end this article by giving some examples of sequences $(g_n)_{n\ge 1}$ for which (\ref{rate}) holds.

Let $g_n(x)= g(\kappa_n(x))$. Consider firstly the case $k=1$, and assume that $g$ is a continuously differentiable function defined in $]-1,1[$.
Clearly, for any $x\in ]0,1[$, $\vert g(x)-g(\kappa_n(x))\vert \le \sup_{|x|\le 1}\vert g^\prime(x)\vert\ n^{-1}$. Consequently, (\ref{rate}) holds. Consider next
the case $k\ge 2$.  Suppose that the function $g$ is differentiable and $\nabla g$ belongs to the Sobolev space $W^{m,p_0}(D)$ consisting of weakly differentiable
functions up to the order $m$, with weak derivatives in $L^{p_0}(D)$. Assume $m>\frac{k}{p_0}$.  By the Sobolev embedding theorem (see \cite{AF03}, page 85, Theorem 4.12),
$\nabla g$ is a continuous function on $D$, and we also get (\ref{rate}).
(See \cite{GM06} for similar remarks for $p_0=2$).

Let $g\in W^{1,p_0}(D)$. Define
$$
g_n(x):= n^k \sum_{j\in
\tilde{I}_n^k}\left(\int_{D_j}g(y)dy\right)\Ind_{D_j}(x), \ x\in D,
$$
where $\tilde{I}_n^k=\{0,\ldots,n-1\}^k$. Applying Equation (7.45) in page 157 of
\cite{GT97}, we have
\begin{align*}
\int_{D_j}\left\vert g(x)-n^k\int_{D_j}g(y)dy\right\vert^{p_0} dx&\le \left(\omega_kn^k \right)^{\left(1-\frac{1}{k}\right)p_0} \left(\frac{\sqrt{k}}{n}\right)^{kp_0}\int_{D_j}\left\vert (\nabla g)(y)\right\vert^{p_0} dy\\
&\le C(k,p_0)n^{-p_0}\left\Vert (\nabla g)\Ind_{D_j}\right\Vert_{L^{p_0}}^{p_0}.
\end{align*}
Since
\begin{equation*}
 \|g-g_n\|_{L^{p_0}(D)}^{p_0}=\sum_{j\in
\tilde{I}_n^k}\int_{D_j}\left\vert g(x)-n^k\int_{D_j}g(y)dy\right\vert^{p_0} dx,
\end{equation*}
we obtain (\ref{rate}).
\begin{rem}
\label{r5.3}
Assume (\ref{rate}). By applying Borel-Cantelli's lemma, as in \cite{GM06}, we can prove the following statements:
\begin{enumerate}
\item With the same assumptions as in Theorem \ref{t5.low}, there exists an a.s. finite random variable
$\xi$ such that
\begin{equation*}
\Vert u-u_n\Vert_{L^2(D)}\le \xi n^{-\nu},
\end{equation*}
a.s., with $\nu\in]0,\frac{1}{2}[$ for $k=1,2$, and $\nu\in]0,\frac{1}{4}[$ for $k=3$.
\item Assume the hypotheses of Theorem \ref{t5.1} and let $\nu$ be as in this theorem. Then, there
 exists an a.s. finite random variable
$\xi$ such that a.s.
\begin{equation*}
\Vert u-u_n^{\varepsilon(n)}\Vert_{L^2(D)}\le \xi n^{-\nu}.
\end{equation*}
\end{enumerate}
\end{rem}
\bigskip

\noindent\textbf{Acknowledgements.}
This paper has been partly done when the first author was visiting the Institute Mittag-Leffler
in Djursholm (Sweden) during a semester devoted to SPDEs.
She would like to thank this institution for the very kind hospitality and support.

The authors would like to thank two anonymous referees for constructive and valuable remarks
that helped to improve the paper.


\begin{thebibliography}{Bi1-81}

\bibitem{AF03} Adams, R.A., Fournier, J.J.F.: \emph{Sobolev
Spaces.} Pure and Applied Mathematics Series. Elsevier, Amsterdam,
2nd Edition, 2003.

\bibitem{adler} Adler, R. J.: \emph{An Introduction to Continuity, Extrema, and Related Topics for General Gaussian Processes}. Institute of Mathematical
Statistics Lecture Notes-Monographs Series, Vol. 12, 1990.

\bibitem{BePa61} Benedek, A., Panzone, R.: \emph{The spaces $L^p$ with mixed
norm.} Duke Math. J. \textbf{28}, 301-324 (1961).

\bibitem{BuPa90} Buckdahn, R., Pardoux, E.: \emph{Monotonicity Methods for White Noise Driven Quasi-Linear
SPDEs.} In: Diffusion Process and Related Problems in Analysis,
Vol. 1 (Evanston, IL, 1989), pp. 219-233. Progress in Probability \textbf{22},
Birkhauser, Boston, MA, 1990.

\bibitem{DM92} Donati-Martin, C.: \emph{Quasi-Linear elliptic stochastic partial differential equation: Markov
property.} Stochastics and Stochastic Reports \textbf{41}, 219-240 (1992).

\bibitem{DN94} Donati-Martin, C., Nualart, D.: \emph{Markov property for elliptic stochastic partial differential
equation.} Stochastics and Stochastic Reports \textbf{46}, 107-115 (1994).

\bibitem{Do84} Doob, J.L.: \emph{Classical Potential Theory and Its Probabilistic
Counterpart.} Grundlehren der Mathematischen Wissenschaften
\textbf{262}, Springer Verlag, New York, 1984.

\bibitem{EH03} Elliott, R.J., van der Hoek, J.: \emph{A general fractional white noise theory and
applications to finance.} Mathematical Finance, Vol. \textbf{13},
No. 2, 301-330 (2003).

\bibitem{GT97} Gilbarg, D., Trudinger, N.S.: \emph{Elliptic Partial Differential Equation of Second
Order.} Grundlehren der Mathematischen Wissenschaften \textbf{224},
Springer Verlag, New York, 1977.

\bibitem{GN96} Gripenberg, G., Norros, I.: \emph{On the prediction of fractional Brownian
motion.} J. Appl. Probab. \textbf{33}, 400-410 (1996).

\bibitem{GM06} Gy\"ongy, I., Mart{\'i}nez, T.: \emph{On numerical solution of stochastic
partial differential equations of elliptic type.} Stochastics: An
International Journal of Probability and Stochastics Processes.
Vol. \textbf{78}, No. 4, 213-231 (2006).

\bibitem{Hu00} Hu, Y.: \emph{A class of stochastic partial differential equations driven by fractional
white noise.} In: Gesztesy, F. et al. ed. Stochastic Processes,
Physics and Geometry: New Interplays II, pp. 317-325. Conference Proceedings,
Canadian Mathematical Society, AMS. Vol. \textbf{29}, 2000.

\bibitem{Hu01} Hu, Y.: \emph{Heat Equations with Fractional White Noise Potentials.} Applied
Mathematics and Optimization \textbf{43}, 221-243 (2001).

\bibitem{HOZ00} Hu, Y., {\O}ksendal, B., Zhang, T.: \emph{Stochastic partial differential equations
driven by multiparameter fractional white noise.} In: Gesztesy, F.
et al. ed. Stochastic Processes, Physics and Geometry: New
Interplays II, pp. 327-337. Conference Proceedings, Canadian
Mathematical Society, AMS. Vol. \textbf{29}, (2000).

\bibitem{HOZ04} Hu, Y., {\O}ksendal, B., Zhang, T.: \emph{General Frational Multiparameter White Noise
Theory and Stochastic Partial Differential Equations.}
Communicatons in partial differential equations. Vol. \textbf{29},
No. 1 and 2, 1-23 (2004).

\bibitem{KRW09} Karniadakis, G. E.,  Rozovskii, B. L., Wan, X.: \emph{Stochastic Finite Element Approximations for
Elliptic Problems with Spatial Gaussian Coefficients.} Preprint,
2009.

\bibitem{K} Krylov, N.V.: \emph{Lectures on Elliptic and Parabolic Equations in H\"older Spaces}. Graduate
Studies in Mathematics, Vol. \textbf{12}, American Mathematical
Society, 1991.

\bibitem{LT91}Ledoux, M., Talagrand, M.: \emph{Probability in Banach
Spaces.} Springer Verlag,\ 1991.

\bibitem{LR08} Lototsky, S.V., Rozovskii, B. L.: \emph{Stochastic Partial Differential Equations Driven by Purely Spatial Noise.}
arXiv:math/0505551v3, 2008.

\bibitem{MSS06} Mart\'inez, T., Sanz-Sol\'e, M.: \emph{A Lattice Scheme for Stochastic Partial Differential Equations
of Elliptic Type in Dimension $d\geq 4$.} Appl. Math. Optim.
\textbf{54}, 343-368 (2006).

\bibitem{mishura} Mishura, Y.: \emph{Stochastic Calculus for Fractional Brownian Motion and Related Processes.}
Lecture Notes in Mathematics 1929. Springer Verlag, 2008.

\bibitem{Nu06} Nualart, D.: \emph{The Malliavin Calculus and Related
Topics.} Probability and its Applications. Springer Verlag, 2nd
Edition, 2006.

\bibitem{PT01} Pipiras V., Taqqu, M. S.: \emph{Integration questions related to fractional Brownian
motion.} Probab. Theory Relat. Fields. \textbf{118}, 251-291 (2000).


\bibitem{R-Y} Revuz, D., Yor, M.: \emph{Continuous martingales and Brownian Motion}. Grundlehren der
Mathematischen Wissenschaften,  \textbf{293}. Third printing of the third edition. Springer Verlag, New York, 2005.

\bibitem{ST94} Samorodnitsky, G., Taqqu, M. S.: \emph{Stable non-Gaussian random
processes.} Chapman and Hall, 1994.

\bibitem{SSV07} Sanz-Sol\'e, M. Vuillermot, P.: \emph{Mild Solutions for a Class
of Fractional SPDEs and Their Sample Paths}. Journal of Evolution Equations, to appear.

\bibitem{Ze90} Zeidler, E.: \emph{Nonlinear Functional Analysis and Its
Applications. II/B. Nonlinear Monotone Operators.}
Springer Verlag, New York, 1990.

\end{thebibliography}
\end{document}